# THE FUNDAMENTALS OF RELATIONS LANGUAGE MATHEMATICS
1. Logic, theory of sets, structures.

V.Ye.Mitroshin
V.N.Karazin Kharkov National University
61077, Kharkov, Svobody Sq., 4

**Abstract**. The fundamentals of formal logic, theory of sets and mathematical structures are narrated in terms of relations language.

The present paper launches a series of three publications under the common title. In its first section relations language [RL], its formal logic [RL-logic], and on this basis RL-theory of sets and mathematical structures are offered. The latter appear as plurality of relations without minimal element, such plurality complying with thereafter defined properties of reflexivity, symmetry and transitiveness. Such pluralities in general case are called "Binary relations filters" or $R_{112}$– filters [1].

Let us instantly invoke a reservation that in RL-logic the notion of "set" presents a classification, but not an object. Null-set also disappears as an *object*, and the notions of *point* and *one-point set* coincide. It has a lot of corollaries. In particular, Zermelo theorem turns out to be unprovable, and presents particular axiom, and the paradoxes, which troubled the founders of the classical trend for so long time, disappear.

In the second section $R_{312}$– filters will be investigated, which induce not only usual topological structures, but also the objects, which do not comply with standard axioms.

In the third section of the paper advanced techniques is used for investigation of various sets of measurable functions and linear densely defined operators.

## Section1. Relations language and its philosophy

Two basic notions – *truth* and *lie* transpierce any sphere of activity of a human being. Philosophic notion of the truth, as *"cognition ideality, which consists in coincidence of conceivable things with reality"*, does not contribute something to mathematics, because the entire "reality" in mathematics is exactly "conceivable". That is why the mathematicians came to an agreement to consider some "state of things" as *evident*, i.e. distinctly observed by mind's eye, without regard to their "materialistic impletion", and consequently as *true*. At this, any statement about not true "state of things" was perceived as *lie*. And over a period of two millenniums development of philosophy and mathematics took place under the motto, that any declarative statement could be treated as true or lying. And though the belief in such "state of things" was every now and then undermined by arising semantic paradoxes, many hoped to find in future such "meta language", where these paradoxes would become just impossible.

In the course of time, obviously, awareness of the fact, that it is difficult to lay the world of thoughts into Procrustian bed of "truth – lie", led not only to emergence of logics, absolutely refusing the principle of the third element exclusion [1,2], but also to emergence of polyvalent "logics" (see review [3]). But the essence of all new trends, to our mind, has remained unchanged – they appear to be successive attempts to narrow down the world of thoughts to certain plurality of algebraic rules, let them be not so elementary, as the code of practice with "true/lying" tables of classic logic.

We suggest a different view on the *language* of mathematics, namely as on the language of *relations*. This view is based on the philosophical idea, which has been sporadically appearing since the times of Eudox, and which in 1902 was clearly stated by Henri Poincare [4]: *"The fact is that it [science] is able to cognize not the essence of a thing, as naive dogmatists think, but only the relations between the things; beyond these relations there is no cognizable reality"*.

One of the most important consequences of this conceptual idea is that it enables to clearly formulate illusory notion of *sense*, which was previously perceived in philosophy as *"inherent logical subject matter of a word, parole, phenomenon,…, the meaning perceived by mind"*.

Sufficiently detailed presentation of the history and philosophy of relations language could be found in [5, 6]. Here we will only try and briefly formulate its basic provisions.

---

[1] The lower indices point at hereinafter defined properties of reflexivity, symmetry and transitiveness in the relations plurality, which induces the structure.

So we accept as guide to action the thought of Poincare, that science is able to perceive not the essence of a thing, but only the relations between the things. Then:

- primary indefinable notions of any theory could be only: "objects" – let us say "points" and "relations" between them, in which they exist. Points and relations are dualistic notions. The point does not exist, if it has no relations with the other points – otherwise we are not able to know something about it; and without points there are no relations, in which they do exist.

Giving name to any object – a point, we imply also those relations, in which it exists with the other names. Some of the relations appear before our mind's eye quite distinctly; about the others we can only guess, and we even have no suspicion about the existence of some. And that is why the "name" is not just some "convenient" combination of letters of the *alphabet*, accepted by us for expression of the names, because the names are *realized* only in the *language*, which expresses the relations between the denominated objects.

We should underline that:

- the "point" is the image of notional relation of equivalency. But those relations, that we perceive, do not discern its "inherent" content. With emergence of new relations the point can appear as a set of "diverse points".

Thus, the following issues appear before our mind's eye: Universal of objects-points and Universal of relations between them – boundless notions. Imaginary horizon of points Universal and relations Universal is vague, not defined, and all the time it moves away from us, as far as we cognize it deeper and deeper. That is why the relations field between the names senses, which constantly expands, forces us to acknowledge, that some portion of our expositions is correct only with "accuracy to …" Together with understanding of new relations between the names, even two previously *equal* meanings of the names, if they are not *identical*, could turn out to be different. Until now, due to historically caused limitation of our knowledge, we just have not noticed the difference existing between the objects: truthfulness of the statement, which is expressed by symbol of *congruence* $x = y$, contrary to the statement, expressed by symbol of *identity* $x \equiv y$, and, being comprehended as different name of "the same", depends on Universal of relations, in which the objects under investigation exist. That is why "variant readings in notions" often occur, and not only in everyday life – even the mathematicians' perception of the relations between "the same objects" may differ a little in some ways. And the statement nonsense for one researcher could turn out for the other one not only sensible, but true. For instance, for logicians statement «$Z$ is the son of childless parents $X$ and $Y$» appears to be absurd, since they *imply* that what is meant here is their biological son. But for another reader this statement is indefinite, in fact it cold refer to an adopted son. In such case this statement will be quite sensible and true. Absurd for the majority of logicians [7,8] statement "I lie" as an answer to the question "Do you lie?", in reality is *sensible* but *indefinite*. For instance, statement "I lie" only on Fridays, or only to my wife has indefinite, hidden sense. But if the point is that "I *always* lie" and it means also right now, then for relations language

- any statement, denying itself, is named absurd.

Often outwardly nonsense statements appear in those cases, when the statement contains the names from different Universals of points. Thus, often cited and outwardly nonsense statement "«Caesar» is a prime number" proves to be quite sensible, if all the names are put in order, and their finite number corresponds to number *N*, and these numbers are put in one-to-one single-valued correspondence to interval of whole numbers $N+1$, then statement "«Caesar» is a prime number" precisely means that name "Caesar" in the adopted sorting of names corresponds to a prime number. Hence,

- in order for the statement to have sense, the "names" in the statement should be either from one Universal of points, or their Universals should be isomorphic.

This is necessary, but by far not sufficient condition of the statements meaningfulness.

We have already many times used the word "sense". To our mind, for any *language*, and, consequently for all the science sub-disciplines in general this notion is principal. Even the most intricate linguists and philosophers can not avoid appealing to "sense", because we write and investigate notions and statements, which we are able to "give sense to" and to transfer this "sense" to the other reader. More than that, even the most genuine formalists in mathematics, for example [9], can not avoid this "sin", and not only on the stage of "tentative explanations" and formulation of axioms (definitions). Just because any symbolic language represents nothing more than abbreviated *connected* record of statements, comprising "tentative" explanations, where the notion of "sense" is sure to be present. But is it possible at all to define the notion of "sense" by means of statements, which should have this sense themselves? In his review of different points of view on the notion of "sense", both linguistic and philosophical ones, A.I.Novikov [10] wrote: *"Sense belongs to the category of those enigmatic phenomena, which are considered to be commonly known, since it plays part both in academic and in ordinary communication. In reality, it is not only deprived of any rigorous generally accepted definition, but even on descriptive level there exists big spread of judgements what it means. Sometimes it is accepted, that the sense belongs to those most general categories, which are not subject for defining and should be perceived as some kind of given entity".*

As a result it turns out that we are standing at very unsettled ground. And the question "Do we really understand each other?" is becoming more and more topical in our technocratic time.

To our mind hopelessness of the situation is conditioned not by perceived, but just by "materialistic" standpoint on the goal of our research – we are still looking for "bricks" of the universe, i.e. the notions, using which we hope to construct the whole World. Again we are searching for the essences, investing them with "qualities", "properties", not taking into account that these notions are relations, but not essences. At this, using the same "bricks" it is possible to construct a prison for the brain or to erect a house of worship for it – everything depends on correlation between the "bricks".

And the "bricks" themselves, their "name that exists" – word root, its spelling is the entire world of relations, which reveals its richness, its many faces in prefixes, endings, contexts....

And though in mathematics the name of an object-point is only a symbol from the alphabet adopted by us, and the math studies small fragments of speculative universe with clear boundaries, as far as possible, both in Universal of points and in Universal of relations between them, it is impossible to guarantee that somewhere in deep subconsciousness some of the relations are only implied (not perceived). And the "sense", embedded in the alphabet and embodied at the stages of "preliminary explanations" and axioms (definitions), is concealed in proved theorem. And it could happen that some of the theorems are deprived of sense at all.

So what is "sensible statement" for relations language, and what is "sense' as such?

- The sense of statement (word) is perceived entirety of relations, in which the named objects exist.

That is why the sense of the statement "point $y$ exists in relation $r$ with point $x$" just consists in the fact that "point $y$ exists in relation $r$ with point $x$". That is why some of the words come to life – they are filled with the sense new for us, while we learn more and more about their relations with the other words. That is why, rereading the old book, sometimes we have the feeling that we have not read it at all. It happens just because the Universals of relations and points, which we perceive, do not have clear boundaries, change all the time – new relations appear, the other, on the contrary, die away, and what was truth for us yesterday, today appears to be lie. Thus,

- we categorize the statements, first of all, as sensible, nonsense (statements deprived of sense) and absurd.

Absurd and nonsense statements usually differ. The first ones are perceived by us, as the statements, which contradict themselves, for example "I always lie"; the second ones are perceived by us as the statements deprived of sense, for example *"If it rains, then the tram"*. But in our perception of sense nonsense statements have no sense. Hence, they present only some part – proper subset of nonsense statements set. Some of nonsense statements, but by no means absurd, could become quite sensible, although in wider field of relations. For instance, the quoted above nonsense statement could also look like that: "*If it rains* heavily in the city, *then the tram* does not go along our street, located in the lowland». And so on and so forth. And that is why

- sensible statements could be true, lying or indefinite.

But speaking about indefinite statements, we have in mind not only the statements, the definition of which could be completed, and they will become either true, or lying, or absurd. For example, the statement «$Z$ is the son of childless parents $X$ and $Y$» becomes absurd, if the definition of word "son" is completed by word "biological"; and it becomes true, if the word "son" definition is completed by word "adopted". We also mean those statements, the definition of which could not be completed theoretically. Such statements often appear as infinite chain of conjunctions $(r) \wedge (s) \wedge (t) \wedge ..$ of statements $(r),(s),(t),...$, each finite chain of which is quite sensible and true. For the sake of example we will present often quoted [11] G.Cantor proof of theorem on nondenumerability of points set of unit segment $E = [0.1]$ of real numbers. It is absolutely simple and very illustrative.

Theorem (Cantor). Points set of unit segment $E = [0.1]$ is nondenumerable.

*Proof.* Let us assume the opposite, that this set countable: $E = \{x_1, x_2, x_3, ...\}$. We select the first segment $s_1$ out of $E$ in such a way, that it should not contain point $x_1$; the second one $s_2 \subset s_1$ in such a way that it should not contain $x_2$; the third one $s_3 \subset s_2$ in such a way that it should not contain $x_3$, and so forth. It could be done every time, as Cantor used to prove, by dividing the segments into three equal parts. In the overlap of the segments, imbedded one into another, i.e. $\bigcap_{i=1} s_i$, according to Cantor and his followers idea, there must for sure be the point, that does not appertain to $E$, and this exactly proves its non-denumerability ⌛

But at this, they leave out of account, that each segment is a subset of denumerable set $E$ according to assumption *ex contrario*, and that is why it also appears to be denumerable set. Hence, overlap $\bigcap_{i=1} s_i$ is not defined, i.e. it is null, because for any point out of $E$ there exist a segment, which does not incorporate it. And though the theorem is true – see Section 3, the above proof is fallible.

The same type infinite chains of conjunctions are not rare in mathematics. In fact, on this basis intuitionists thought [1] that the principle of the third element exclusion can not be a law of logic, overlooking that a statement could turn out to be indefinite.

In this connection we would like to draw your attention to the fact, that in mathematical theory existence of the statement in a form of indefinite infinite conjunctions chain in substance presents Godel theorem about existence in this theory of the statements, which can not be either proved, or invalidated. We will speak just a little later about formal facet of the logic of relations language – about formalization of the notion of "sense". And we will finish the discussion by general characteristic.

- For relations language declarative statements could be sensible and nonsense. Nonsense statements in the capacity of proper subset contain absurd statements. Sensible statements could be true, lying and indefinite. Under expansion of relations Universal a) some of the indefinite statements could turn either true, or lying, or absurd, and b) some of nonsense statements, excluding absurd ones, could turn sensible.

Therefore, logic of the relations language will turn out to be "four-dimensional", but not two-dimensional "yes/no" of classical calculus. At this, it will be as if "alive", often changing, depending on mathematical theory, being developed.

## Section 2. RL-logic

**Designations**

In RL-logic (Relations Language Logic) points and a relationship in which they exist, are marked by the Latin lowercase letters, but in different styles.

Symbol $\exists$, standing in front of point $x$ (relation $r$), means: *there exists $x$* (respectively $r$); denial of existence of the point (relation) is designated as $\nexists$. Record $\exists!...$ means *singularity* of the object, the existence of which is affirmed.

Statement $\forall x...$ implies: *for each point $x$ ...,* and each semi-column ... **:** ... – *thus, (of the kind), such, that...*

Symbols $\wedge$ and $\vee$ designate conjunctions *and,* and respectively *or.*

Unless otherwise noted, the equality symbol "=" is understood as the definition of the left part of the statement – the part, that stands ahead of symbol "="; through the right part – the part, that follows symbol "=" is defined. Relation of identity, "≡", unlike equality, is not changed at any extension of the Universals.

In this section, for simplicity of speech, we will use the standard symbols of set theory. The numbering of the statements, theorems and examples is consecutive. The first digit means the section number; the second is the sequential number. End of examples, deductions and comments will be represented by the graphic symbol of the sand glass – ⌛.

**Basic notions**

Further on, meanings of the names – objects $x, y, z,..$, their certain "field" $F = \{x, y, z,...\}$, will be called by us not only as *points*, but as *primitives*, by that pointing out non-definability of these notions, except by means of the "field" of $n$-ary relations $R \equiv \{r, s, t,...\}$, in which these points exist among themselves. We assume relations field $R$ as *defined* and *non-contradictory*. Speaking about determinacy of $R$, we imply, that each relation $r$ out of $R$ is determined in the whole field of primitives. Speaking about non-contradiction of $R_=$, we imply that in $R_=$ there is no pair of relations, which are incompatible with each other. We will formalize these notions thereafter. We will call dyad $\{F,R\}$ *the object of logic*.

Speaking that *point y* exists *in relation* $r$ *with point x*, we will write $(x[r]y)_{Tr}$, stressing by lower symbol truthfulness of the statement (derivative from **tr**uth).

Record $(x[r]\{y,z,..\})_{Tr}$ means that that plurality of points $\{y,z,..\}$ (the sequence order of which and their number $n$ are fixed by relation $r$) is in relation $r$ with point $x$ ($n$-ary relation). If we speak about lie or equivocation of the statement, we write $(x[r]\{y,z,..\})_{Li}$ and $(x[r]\{y,z,..\})_{Pr}$, respectively (derivatives from **li**e and **pr**imitive). When there is no need to indicate particular primitive statement, then, speaking about *statement*, we will just frame relation $r$, inducing it, by round or square brackets, i.e. we will write it in the form of $(r)$ or $[r]$.

All the points $\{y,z,..\}$ out of $F$, for which $(x[r]\{y,z,..\})_{Tr}$, are called *image* of $r$ in point $x$, that is recorded as $r[x]=\{\{y,z,..\}:(x[r]\{y,z,..\})_{Tr}\}$.

All the points $x$ out of $F$, for which image $r[x]$ is not "null", and which is designated by symbol "$\phi$", i.e. $I_r=\{x:r[x]\neq\phi\}$, is named truthfulness domain of $r$. And since each relation out of $R$ *is defined* for the whole $F$, then addition to $I_r$ – all the points out of $F\setminus I_r$ is entitled as the domain of *lie*. If $I_r=\phi$ is null, then $r$ is called relation undefined in $F$.

Grammatically properly constructed statement $(r)$ can have sense – be sensible $(r)_{Se}$, but it can have no sense – be nonsense $(r)_{No}$ (derivatives from **se**nse and **no**nsense respectively). Sensible statement can be true $(r)_{Tr}$, lying $(r)_{Li}$, or undefined $(r)_{Pr}$. We will name statements $(r)$ and $(s)$ *equipollent* and write $(r)\square(s)$, if they express the same sense. Consequently, $I_r=I_s$.

*Example*. To illustrate efficiency of the adopted symbols, let us look at the classical paradox of the "liar" – "I always lie". Now we formalize this statement. We will eliminate, first of all, "I" – we do not care who particularly is speaking. The word "always" is equipollent to the word "any". And the statement "I always lie" will be written as $[(\forall r)_{Li}]_{Tr}$ – any statement is lying. But for relations language, to affirm truthfulness of lying state of things or lie of the true state of things, for example, $x\neq x$, is just an absurd.

On the other hand, *sensible* statement "I lie" as an answer the question "Do you lie?" is indefinite, but not an absurd statement, as some textbooks assume: $\exists r:(r)_{Li}$ – «something», «sometimes» and etc., happens to be lying ⌛

Thus, we distinguish statements, first of all, as sensible and nonsense. But the notion of sensibility, previously formulated, should be formalized.

Let $R$ be plurality of non-contradictory relations, perceived by our mind and defined by us, in which primitives out of $F$ exist. These relations are just integral to the definitions of the notions of the mathematical theory, which we intend to develop. That is why these relations are *a priori* sensible, and this is recorded by us as the first axiom of relations language.

**2.1**. Axiom. Any relation of $r$ out of $R$ is sensible.

And since any relation $r$ out of $R$ has its own *not null* truthfulness field $I_r$, then, notion of *determinacy* of relation $r$ consists just in this, and as addition in $F$, there exists lie field $F\setminus I_r$, so we come to a conclusion, that

**2.2**. Statement. For each $r$ out of $R$ the principle of the third option exclusion is correct: $(r)_{Tr}\vee(r)_{Li}$.

Basing on a priori relations sensibility of field $R$, we will be interested to find out the conditions of *sensibility* of these relations combinations, which rely on conjunctions $\wedge,\vee$, implication symbol "$\rightarrow$" and equivalency symbol "$\leftrightarrow$". This is the first task of assertions logic.

Which way relation $r\wedge s$ should be interpreted?

**2.3**. Definition. Relation $r\wedge s$, which is entitled as conjunction, is understood as plurality of statements
type: $(x[r\wedge s]\{y,z,..\})_\lambda=(x[r]\{y,z,..\})_\mu\wedge(x[s]\{y,z,..\})_\nu$

Here symbols $\lambda,\mu,\nu$ signify one of the symbols $Tr, Li, Pr, Ab$. Let us look at it more in details. Each of the statements, comprising the right part, is induced by sensible relations $r$ and $s$. That is why indices $\mu, \nu$ can acquire the only valuations Tr and Li. In set $I_\wedge = I_s \cap I_r$ of conjunction truthfulness, if it is not null, $\lambda = \mu = \nu = \text{Tr}$. In the other points the statement is lying. If $I_\wedge$ is null, then, disregarding whatsoever set $F \backslash (I_s \cup I_r)$ could be, we will always deal with either indefinite, or absurd statement.

**2.4**. Example. Hereinafter in all the examples $F = N$ will mean a set of natural numbers. We will mean by record $m/k$, that $m \in N$ is divided by $k \in N$ without remainder.

**1**. Let us consider two binary relations:

$$(m[\text{r}]n) = [(n = m) \wedge (m/2)], \quad (m[\text{s}]n) = [(n = m) \wedge (m/3)].$$

These are two sensible statements. The first one defines the numbers, which are divisible by 2, the second one – those, which are divisible by 3. Truthfulness domains of these statements do not coincide, though they overlap. Now let us consider compound statement $(\text{r}) \wedge (\text{s})$:

$$(m[\text{r} \wedge \text{s}]n) = ([n = m] \wedge [(m/2) \wedge (m/3)]).$$

It is also sensible, it defines the numbers, which are divisible by 2 and by 3.

**2**. Let us have a look at two binary relations

$$(m[\text{p}]n) = ([n = m] \wedge [(-1)^n < 0]), \quad (m[\text{q}]n) = ([n = m] \wedge [(-1)^n > 0]).$$

These are two sensible statements. The first one defines odd numbers, and the second one – even numbers. Their truthfulness domains do not coincide and do not crosscut. And here comes compound statement

$$(m[\text{p} \wedge \text{q}]n) = [(n = m) \wedge ([(-1)^n < 0] \wedge [(-1)^n > 0])],$$

which is not only not defined (has null truthfulness domain), but at the same time it is absurd, since one of the parts of the sentence refuses the second one. And that is what we mean by absurdity. We have to point out, that in standard calculus for any "formulas" p and q expression for $p \wedge q$ is also considered to be "formula". For relations language, as we see, it is not always the case ⌛

That is why we will assume the following in the capacity of the second axiom of relations language logic:

**2.5**. Axiom. Relation $\text{r} \wedge \text{s}$ is sensible, if and only if $I_r \cap I_s \neq \phi$.

Truss $\wedge$ makes it possible to give definition of *non-contradiction* of field R: $\exists \text{r,s}:[(\text{r})_{\text{Tr}} \wedge (\text{s})_{\text{Tr}}]_{\text{Li}}$.

**2.6**. Definition. Relation $\text{r} \vee \text{s}$, which is called *disjunction*, is understood as plurality of the statements type: $(x[r \vee s]\{y,z,..\})_\lambda = (x[r]\{y,z,..\})_\mu \vee (x[s]\{y,z,..\})_\nu$

**2.7**. Example. Relations r,s and p,q are the same as in 2.4. Both statement

 **1**. $(m[\text{r} \vee \text{s}]n) = ([n = m] \wedge [(m/2) \vee (m/3)])$,

and statement

 **2**. $(m[\text{p} \vee \text{q}]n) = [(n = m) \wedge ([(-1)^n < 0] \vee [(-1)^n > 0])]$

have univalent sense. It seems to be impossible to make nonsense statement with conjunction *"or"*. But why then both statement [12]

 **3**. « $2 \times 2 = 4$ », *or* "New York is a big city",

and statement

 **4**. « $2 \times 2 = 5$ », *or* "New York is a big city",

which are *absurd* for common sense, prove to be *true* in classic logic?

Let us formalize these statements. In set of natural numbers statements $2 \times 2 = 4$ and $2 \times 2 = 5$ are induced by relation $t(k)$ type: $(m[t(k)]n) = (n \times m = k)$, where $k = 4,5$. Then, $I_{t(4)} = \{1,2,4\}$, and $I_{t(5)} = \{1,5\}$.

Let us formalize statement "New York is a big city". By word "big" or "small" we will mean, for instance, the area $S_{name}$ of the city (*name*) up to the city boundaries. The number of the cities on Earth is finite; let it be $N$. That is why isomorphism f exists between subset of natural numbers $\mathbf{Z} = \{k : k \leq N\}$, and set of the cities, arranged in the order of increasing of the areas they occupy, i.e. $k < l$ result in $S_k < S_l$. Let us split $\mathbf{Z}$ into three non-overlapping sequential parts, which will be designated as $\mathbf{Z_1}$, $\mathbf{Z_2}$ and $\mathbf{Z_3}$. If $k \in \mathbf{Z_1}$, then the cities complying with these values of $k$ will be called "small"; the cities with $k \in \mathbf{Z_2}$ will be called "medium". And the cities with $k \in \mathbf{Z_3}$ will be called "big". Let us assume,

$$(m[v(\ell)]n) = [(n = m) \wedge (n \in \mathbf{Z}_\ell)], \text{ c } I_{v(\ell)} = \mathbf{Z}_\ell.$$

Relation $t(k) \vee v(\ell)$ is quite sensible in $\mathbf{Z} = \{k : k \leq N\}$, with non-null truthfulness domain. But, in view of definition 2.6, example 2.7.3 will be recorded like this: $(2[t(4)]2) \vee (2[v(3)]2)$. But, obviously, New York holds not the second number in the cities classification with respect to their areas, and statement $(2[v(3)]2)$ turns to be lying, but not true, as it is affirmed in statement 2.7.3. That is why the second part of sentence 2.7.3 will be recorded like this: $((2[v(3)]2)_{Li})_{Tr}$. But, as we mentioned before, the statements of this kind are absurd. The same way, in example 2.7.4 we have $(2[t(5)]2) \vee (2[v(3)]2)$. But its first part is lying, the second is absurd.

Here the reason is revealed, due to which statements 2.7.3-4 are really absurd, though they are true in classic logic. Example 2.7.3 is built up in the following way. Its first part «$2 \times 2 = 4$» is taken from domain $I_{t(4)}$ of truthfulness $t(4)$, and the second part – "New York is a big city" is taken from domain $I_{v(3)}$ of truthfulness $v(3)$ with $I_{t(4)} \cap I_{v(3)} = \phi$. And these two randomly selected statements were united by conjunction "or" in contravention with definition 2.6 ⌛

The possibility of existence of statements, contradicting common sense, the statements appearing due to conjunction *"or"* and *implication* (see below), at its early stage really perplexed the minds of logicians to big extent [12,13]. But, without solving this problem, they just "turned a blind eye to it". We are solving this problem having defined all the elements of the truss in one field of primitives (or isomorphic fields), and by calculating according to definition 2.6. That is why we accept the following as the third axiom of relations language:

**2.8.** Axiom. For any $r, s$ out of R relation $r \vee s$ is sensible.

**2.9.** Definition. Relation $r \to s$, which is read "if (r), then (s)" and called *implication*, is understood as plurality of statements type:

$$(x[r \to s]\{y,z,..\})_\lambda = (x[r]\{y,z,..\})_\mu \to (x[s]\{y,z,..\})_\nu.$$

Due to determinacy of r and s indices $\nu, \mu$ acquire only valuations Tr or Li. And since we are interested only in true implications, we have to find out, under which conditions index $\lambda$ acquires valuation Tr.

We should mention, that, if $(r)_{Tr} \to (s)_{Tr}$, then $I_r$ should coincide with $I_s$; otherwise in some subset of true values of (r) statement (s) can turn out to be lying, and in some subset of lying values of (r) statement (s) will remain true. That is why condition $I_r = I_s$ is *necessary* condition for the implication to be true. Analogous considerations show, that, if $(r)_{Li} \to (s)_{Li}$ is true implication, then $F \setminus I_r = F \setminus I_s$ should take place, which is again tantamount to $I_r = I_s$. Whenever $(r)_{Li} \to (s)_{Tr}$, we get $I_s = F \setminus I_r$.

The above three types of statements are referred to true implications.

**2.10.** Example. $F = N$. Let
$(m[r]n) = [(n = m) \wedge (m/2)]$, $(m[p]n) = ([n = m] \wedge [(-1)^n > 0])$, $(m[q]n) = ([n \geq m] \wedge [(-1)^m > 0])$.

**1.** Implication $r \to p$ is true with $I_r = I_p$, both in set of even numbers, then $(r)_{Tr} \to (p)_{Tr}$, and in set of odd numbers, then $(r)_{Li} \to (p)_{Li}$

**2.** Implication $r \to q$ c $I_r = I_q$ is also true: a) if $n = m$, then $(r)_{Tr} \to (q)_{Tr}$; b) if $n > m$, then $(r)_{Li} \to (q)_{Tr}$; c) if $n < m$, then $(r)_{Li} \to (q)_{Li}$ ⌛

When truthfulness of (r) results in lie of (s), we get $I_r = F \setminus I_s$. In standard calculus statements $(r)_{Tr} \to (s)_{Li}$ are a priori referred to lying implications, which is not always the case.

**2.11.** Example. $F = N$. Let
$(m[r]n) = [(n = m) \wedge (m/2)]$ и $(m[q]n) = ([n = m] \wedge [(-1)^n < 0])$.

Let us consider implication $q \to r$, i.e. $([n=m] \wedge [(-1)^n < 0])_\mu \to [(n=m) \wedge (m/2)]_\nu$. In subset of odd numbers implication $(q)_{Tr} \to (r)_{Li}$ is true, since it affirms that, if "$n$ *is an odd number*" is true, then "$n$ *is divisible by* 2 *without remainder*" is lying. In subset of even numbers implication $(r)_{Li} \to (q)_{Tr}$ is also true, affirming, that, if "*n is an odd number*" is lying, then "*n is divisible by* 2 *without remainder*" is true ⌛

We see that there are no grounds, except for apriori ones, to consider implication $(r)_{Tr} \to (s)_{Li}$ undoubtedly lying – everything depends on correlations between $I_r$ and $I_s$. That is why we adopt the following as the fourth axiom:

**2.12.** Axiom. Relation $r \to s$ is sensible, if and only if $I_r = I_s$, or $I_s = F \setminus I_r$ ($I_r = F \setminus I_s$).

The relations, which do not comply with the above axioms, have no sense for relations language, i.e. they are nonsense statements.

**2.13.** Example. $F = N$. Let us look at the relations:

**1**. $(m[p]n) = ([n=m] \wedge (-1)^n > 0])$, $(m[t(4)]n) = (n \times m = 4)$.

Can any kind of implication do exist between them? In standard calculus it can without any doubt. In relation logic it can not, since $I_p = \{2,4,6,...\}$ and $I_t = \{1,2,4\}$, and axiom 3.12 is not fulfilled. That is why any implication between these relations has no sense.

**2**. In standard calculus type *"if* 2 > 3*, then* 138 *is a prime number"*, is true assertion. We formalize this statement in $F = N$. Subset of prime numbers will be designated as Dr. Statements type 2 > 3 are brought to existence by relation $(m[w]n) = (m > n)$, with image $w[m] = \{1,2,...,m-1\}$ and truthfulness domain equal to $I_w = N \setminus \{1\}$. Statement "138 is a prime number" is generated by relation $(m[d]n) = [(n=m) \wedge (n \in D)]$, with truthfulness domain $I_d = Dr$. It is clear, that axiom 2.12 is not fulfilled, and consequently *"if* 2 > 3*, then* 138 *is a prime number"* is nonsense statement.

**3**. In standard calculus statements type "if 2 > 3, then New York is a big city" are true assertions, but they are absurd in RL-logic. In order to make sure in this, it is necessary to model reflections from 2.7. ⌛

**2.14.** Definition. Relations r,s are called equivalent and recorded $r \leftrightarrow s$, if and only if $(r \to s)_{Tr} \wedge (s \to r)_{Tr}$, which is understood as plurality of statements type:

$$(x[r \leftrightarrow s]\{y,z,..\})_\lambda = [(x[r]\{y,z,..\})_\mu \to (x[s]\{y,z,..\})_\mu] \wedge [(x[s]\{y,z,..\})_\mu \to (x[r]\{y,z,..\})_\mu]$$

Verbal proof, analogous to that, which we pursued when investigating one-sided implication, show, that index $\lambda = Tr$ in case and only in case, when $I_r = I_s$, both for $\mu = Tr$, and for $\mu = Li$. Then we say, that statements (r) and (s) *are equipollent*.

2.10.1 gives us an example of equivalent relations and equipollent statements; 2.10.2 and 2.13 are the examples of inequivalent relations.

We pay attention, that identity sign "≡" is not tantamount equivalency sign "↔". Relations $(m[r]n) = [(n=m) \wedge (m/2)]$ and $(m[p]n) = ([n=m] \wedge [(-1)^n > 0])$ are equivalent, but *a priori* are not identical.

To affirm lie of true statement and truthfulness of lying one is absurd. It is possible to convince oneself in this by direct calculation.

**2.15.** Statement. Statements $((r)_{Tr})_{Li}$ and $((r)_{Li})_{Tr}$ have no sense.

*Proof*. Let us look at the first statement. Since $((r)_{Tr})_{Li} \equiv (r)_{Tr} \wedge (r)_{Li}$, then, assuming that $(s)_{Tr} = (r)_{Li}$ with $I_s = F \setminus I_r$, we get $I_r \cap I_s = I_r \cap F \setminus I_r = \phi$, which contradicts axiom 2.5. And since $(r)_{Tr} \wedge (r)_{Li} \equiv (r)_{Li} \wedge (r)_{Tr}$, we get the proof for the second part of the statement ⌛

**Relations field structure**

Now we will discuss the issue about structure of field **R**. First of all, we should note, that the principle of the third option exclusion is fulfilled for any point r of field **R**. And, consequently, it should be fulfilled also for any element equivalent to r. That is why, in future, by r we will imply its *equivalence class*. Hereinafter for any *given* conjunctions chain $r \wedge s \wedge t \wedge ...$ the principle of the third option exclusion will also be fulfilled. This principle is also applicable to $r \vee s \vee t...$ It follows that such kind of relations can be the elements of field **R**, and let us add them to it. But implication $r \rightarrow s$ is external relation over field **R**. It can also be absurd. That is why such kind of relations can not be the elements of field **R**.

In field **R** there can be relations, represented just in a form of infinite chain of conjunctions. If $I_z = \phi$, then it means lie of statement (z) in the whole field of primitives. And these are indefinite statements, which are not available in field **R** by definition. If $I_z \neq \phi$, then either $(z)_{Tr}$, or $(z)_{Li}$ take place. That is why for field **R** the principle of the third option exclusion is fulfilled without any limits. But field **R** can not be considered complete, since there could be a chain of conjunctions, the limit of which has null field of truthfulness.

Extension of field **R**, for example, perception of new relations can also lead to extension of plurality of sensible statements – some of previously nonsense statements will turn sensible. It gives all the grounds to assume, that logic will never be able to become completed theory and will forever remain a "stepchild" of ever developing mathematics.

### The laws of logic

In standard calculus (including intuitionists) the following assertions are considered to be *the laws of logic*, arising from definition of the respective trusses, in particular:

C1. $\forall r,s\ [(r \wedge s) \rightarrow (r)]$; $\forall r,s\ [(r \wedge s) \rightarrow (s)]$.

C2. $\forall r,s\ [(r) \rightarrow (s \vee r)]$; $\forall r,s\ [(s) \rightarrow (s \vee r)]$.

C3. $(r) \rightarrow [(\forall s) \rightarrow (r)]$.

C4. $(r) \rightarrow [(r) \rightarrow (\forall s)]$.

C5. $\forall r,s\ [(r) \rightarrow (s)] \vee [(s) \rightarrow (r)]$.

Let us consider assertions C1. In compliance with 2.5, sensibility condition for the left part of statement $(r \wedge s) \rightarrow (r)$ is $I_r \cap I_s \neq \phi$. Then, in compliance with 2.12, sensibility condition for implication between the assertions will be condition $I_r \cap I_s = I_r$, which in reliance on the second statement C1 produces the following: $\forall r,s\ I_r = I_s$. But for random relations this condition is unsatisfiable. That is why for relations language C1 can not act as a logic law.

Let us consider assertions C2, which seem to be quite compelling. The area of true values of $s \vee r$ will be designated as $I_\vee$. Then, the condition for truthfulness of implication $r \rightarrow s \vee r$, in accordance with 2.12, can only be $I_r = I_\vee$. The other possibilities are excluded, since in the right part of the implication the same relation as in the left part is present. It instantly results in $I_r \supseteq I_s$. Analogous reasoning for the second implication C2 produces $I_r \subseteq I_s$, and that is why the binding condition of C2 truthfulness is $I_r = I_s$. But it is impossible for random relations $r, s$.

Let us consider assertion C3, which sounds pretty strange in colloquial language: *true assertion proceeds from any assertion*. Then, on account of axiom 2.12, we find the following conditions for implication truthfulness: $(\forall s)_{Tr,Li} \rightarrow (r)_{Tr}$: either $\forall s\ I_r = I_s$, or $\forall s\ I_r = F \setminus I_s$, which is impossible for random $s$. Hence, C3 also can not be regarded as a law of logic.

Analogous verbal proofs show, that both C4, and C5 contradict axiom 2.12. That is why we can adopt only the following in the quality of the *laws* of logic of relations language (generality quantifier is not specified, but implied):

R1. of identity: $(r) \equiv (r)$.

R2. of the third option exclusion: $(r)_{Tr} \vee (r)_{Li}$.

R3. of absurdity: $([(r)_{Tr}]_{Li})_{Ab}$, $([(r)_{Li}]_{Tr})_{Ab}$.

R4. of double negation: $((r)_{Li})_{Li} \equiv (r)_{Tr}$.

R5. of contraposition: $(r)_{Tr} \rightarrow (s)_{Tr} \equiv (s)_{Li} \rightarrow (r)_{Li}$.

R6. of Morgan: $(r \wedge s)_{Li} \equiv (r)_{Li} \vee (s)_{Li}$, $(r \vee s)_{Li} \equiv (r)_{Li} \wedge (s)_{Li}$.

R7. of transitivity $(r \to s)_{Tr} \wedge (s \to t)_{Tr} \to (r \to t)_{Tr}$.

## About deductions

In contrast to deductions [14,15] of classical calculus, we do not distinguish between the notions of "implication" (implication – consequence, conclusion) and "sequent" (sequence – succession), though the difference in *conceptual sense* of expressions "if $(r)$, then $(s)$" and, respectively, "$(s)$ arises from $(r)$" formally exists due to ill-defined "if". But since each *true* implication in relations language always begins more specifically – either $(r)_{Tr}$, or $(r)_{Li}$, then this conceptual difference disappears.

In reference to the notion of "formulas deducibility" according to the rule *modus ponens* and *synthesis*, which are systematically used in predicates calculus, then for relations language these are nothing more than tautologies. For example, *modus ponens* is recorded this way: $(r) \wedge [(r) \to (s)] \to (s)$.

What is *theorem* (any) for relations language? It is nothing more than new plurality of relations in F, in sensibility of which we would like to make sure. That is why the deduction itself of the truthfulness (falsehood) of some kind of statements – theorem (hypothesis) in relations language always represents some transitive-reflexive relation over field R, generated by implications plurality. Naturally, this chain of logisms is based on definitions of logic trusses, generality and existential quantifiers, as well as logical laws. And, consequently

- each proved theorem and its conclusions, in essence, present existing, but previously not perceived by us relations in F.

By this, the development of the *theory of proofs* is equally matched to development of the theory of *partially ordered structures* and does not have solution perspective – there will be nothing to say except for general phrases. The author hopes to discuss the elements of theory of proofs in special study. But we would like to highlight some issues right now.

The proof of theorem $(t)$ truthfulness may require proof of intermediate statement, we will say $(z)$. In its turn, it may require consideration of infinite chain of conjunctions: $(z) = (r) \wedge (s) \wedge ...$ If it turns out that $I_z = \phi$, then such statement is not defined. But any implication $(z)_{Pr} \to (t)_{Tr}$, by reason of axiom 2.12, is absurd, and that is why theorem $(t)_{Tr}$ is unprovable. In Section 1 we made an example of fallible reasoning of this kind, though there are much more of them.

The other important thing is that in proofs "by contradiction" we assume lie of some statement $(t)_{Li}$. Later one strives to demonstrate, that it ends in conclusion about lie of a fortiori truthful statement $(s)_{Tr}$, i.e. results in absurd: $(t)_{Li} \to ... \to ((s)_{Tr})_{Li}$. From here the conclusion is drawn, that the statement is true – $(t)_{Tr}$. But in the process of deduction, at some intermediate stage truss $(p)_{Tr} \vee (q)_{Tr}$ may appear with $I_p \cup I_q = F$. Here the implication chain may terminate. Hence $(t)_{Tr}$ is unprovable.

## Conclusions

Any symbolic language is nothing more than concise presentation of statements, *bound* by the laws of statements logic. But if different people even to small extent, but differently percept the sense of these statements, then, about what "strict language" can we talk, if the sense is not defined? And in this respect Gilbert program was doomed to failure from the very beginning. But the main thing lies in the fact, that the notion of "true" statement is only one of the three hypostasis of the notion of "sensible" statement.

Formal theory of sense, narrated above, is based only on a priori sensibility of the "primary" relations out of field **R**. Why we so "easily" agreed to that? Just because affirmative statement *"point y exists in relation r with point x"* exactly means its sensibility. And then, the only thing remaining for us was to investigate the conditions of sensibility of relations received out of field **R** by means of trusses $\wedge, \vee, \to, \leftrightarrow$. And these things are revealed in the same field of primitives **F** (or isomorphic fields).

It instantly leads to the situation, when we are not able to unify operations $\wedge, \vee, \to, \leftrightarrow$ in a form of "truthfulness tables" for random sensible relations r and s. And moreover, we can not unify the operation of negation $\bar{r}$, which initially was not introduced, because for relations language it is indefinite statement: in fact *"not r"*, means *"some other relation"*. As a result, it is impossible to narrow down formal theory of sense to any kind of plurality of "algebraic" rules due to axioms 2.5,12. Logic is converted into "topological" matter with "algebraic" suspension. In this connection

- building-up of any kind of mathematical theory after defining the field of primitives and relations, in which they exist, requires construction of the accompanying logic of assertions.

## Section 3. Theory of sets in RL-logic

According to George Cantor "The set is coufounding of particular different objects of our intuition or of our brainwork into one integral. These objects are called the points of obtained set". But after paradoxes of the theory of sets were found, a great deal of researchers of the end of the XIX and the beginning of the XX century started looking for more satisfactory definition of the theory of sets (see, for instance[16]). And in the first quarter of the XX century the standpoint was formed, which is also being "professed" nowadays [16,17], that the notion of the "set" should be defined through undefinable *pertaining relation*: "$X$ is a set in the case and only in the case, if $\exists Y : X \in Y$". But this is not definition of the set, but additional condition in the definition, which was never formulated, and which enables to avoid serious paradoxes, and to preserve the property of the set to be a member of "something else". Maybe, that is why, upon sufferance from all sides, A.V. Arkhangelskiy has written [18]: *"The experience of present-day mathematics and analysis of its foundations show, that the sets are that basic elementary material, out of which all the main mathematical objects are constructed"*. And that: *"...there exist all grounds to consider the idea of sets one of the most crucial and the most primary forms of thinking"*.

Intuitive-descriptive approach to the notion of the set is preserved in investigations of the past decades – in the alternative theory of sets developed by P.Vopenko [19].

But not only for the researchers, but also for the historians of science the other point of view on the notion of the "set" remained unnoticed, and which was formulated by Poincare at the beginning of the XX century [4]: *"To define a set always means to make classification, to isolate the subjects, which belong to this set, from those which do not participate in it"*. In other words, for Poincare the notion of the "set" is the classification of Universal of points by some relation.

Such opinion about the notion of the "set" naturally followed Poincare's philosophical view on cognitive essence of the science. Let us remind that: *"The fact is that it [science] is able to cognize not the essence of a thing, as naive dogmatists think, but only the relations between the things; beyond these relations there is no cognizable reality"*.

Many Poincare's coevals, especially G.Weyl [20], brought out superficially close opinions on the notion of the set *"Nobody can describe infinite set differently than, having specified the properties, intrinsic for the elements of this set; ... Conceptualization of infinite set as a kind of plurality, composed by means of infinity of separate arbitrary selection actions, and contemplated later by our consciousness as something integral, is devoid of sense;"inexhaustibility" is embedded just in the essence of infinity"*. Similar thing was said also by Bourbaki [21]: *"The set is formed by the elements, capable of having some properties and exist among themselves or with the elements of the other sets in some kind of relations"*. But, according to Bourbaki, the "set" is always a "member" of the other set.

In theory of sets *relation* is some subset of Cartesian product of sets; and in relations language it is just *set*, which represents *classification*, or, sometimes we will say, Universal points *filtration* by some relation. But it turned out, that even such simple shift of emphasis results in substantial differences with the results of theory of sets. Fortunately, practically all these differences lie in "ultramundane" domains of the mathematics.

**Basic notions**[1)]

Let us fix certain relation $i_\alpha$.

**3.1 Definition.** All the points of Universal $F$, for which $(x[i_\alpha]x)_{Tr}$, are called plurality (set) $X_\alpha$ of points in terms of $i_\alpha$, and recorded as: $X_\alpha = \{x : (x[i_\alpha]x)_{Tr}\}$.

---

[1)] Giving an account of the fundamentals of RL-theory of sets, we, in order not to refer the reader to the respective literature all the time, in our comments and footnotes we will reveal interrelation with axiomatic approach by Zermelo-Fraenkel [13] in theory of sets (ZF-theory).

It means, that we just indicate those points $x$ of Universal F, for which statement $(x[i_\alpha]x)$ is true. They are these singled out points, which we call the set [1].

Pertaining relation "$\in$" becomes dependent notion: truthfulness (lie) of statement $[x \in X_\alpha]$ is determined by truthfulness (lie) of statement $(x[i_\alpha]x)$, i.e. $[x \in X_\alpha]_{Tr,Li} = (x[i_\alpha]x)_{Tr,Li}$

Relation $i_\alpha$, reflexive properties of which determine set $X_\alpha$, will be called *identifier* of this set.

Identifier $i_\alpha$ is *defined* only in the case, if on the visible horizon of Universal of points $\exists x : (x[i_\alpha]x)_{Tr}$. Otherwise we can not know anything about it. That is why for relations language it is impossible to define "null-set" by statement $X_\alpha = \{\overline{\exists} x : (x[i_\alpha]x)_{Tr}\}$, it just means, that identifier is not defined.

**3.2**. Commentary. The concept of set, "which does not contain any elements", plays in ZF-theory fundamental but not auxiliary part, and it is defined in the following way: $\phi = \{x : x \neq x\}$ ⌛

But for relations language, as we said before, the relations between "something" and "nothing" do not exist, and consequently the notion "null-set" does not exist too. But this is philosophy. However, in RL-logic the problem is solved in a pretty simple way: predicate $(x \neq x)$ is given by $[(x \equiv x)_{Tr}]_{Li}$, and consequently it is not absurd. Therefore,

**3.3**. Statement. In RL-theory of sets, null-set does not exist as an object.

Identifier $i_\alpha$ can also be such, that plurality $X_\alpha = \{\exists ! x : (x[i_\alpha]x)_{Tr}\} \equiv \{x_\alpha\}$ will incorporate one point. But in RL-logic the notions of "point" and "single-point set" coincide, since they express the same statement: $\exists ! x : (x[i_\alpha]x)_{Tr}$.

**3.4**. Comment. In this connection usual for mathematicians objects type $\{x, \{x\}, \{\{x\}\}, ..\}$ loose their sense, as well as *axiom of infinity* [2] of ZF-theory ⌛

Equality of two sets in ZF-theory is defined by axiom of *equal volume* (extensionality) [3]. For relations language it is done through RL-equivalency.

**3.5** Comment. Sets $X_\alpha$ and $X_\beta$ are called *equal* and we write $X_\alpha = X_\beta$, if and only if
$\forall x \ (x[i_\alpha]x)_{Tr,Li} \leftrightarrow (x[i_\beta]x)_{Tr,Li}$.

In other words: the sets are equal, if their identifiers are RL-equivalent.
For the sake of completeness let us consider now sets manipulations.

**3.6** Definition.
1. Set $X_\alpha \cup X_\beta = \{x : (x[i_\alpha]x)_{Tr} \vee (x[i_\beta]x)_{Tr}\}$ is called unification of sets $X_\alpha$ and $X_\beta$.
2. Set of points $X_\beta$, which do not belong to $X_\alpha$, i.e. $X_\beta \setminus X_\alpha = \{x : (x[i_\beta]x)_{Tr} \wedge (x[i_\alpha]x)_{Li}\}$ is called complement of set $X_\alpha$ up to set $X_\beta$.
3. Set $X_\alpha \cap X_\beta = \{x : (x[i_\alpha]x)_{Tr} \wedge (x[i_\beta]x)_{Tr}\}$ is called overlap of sets $X_\alpha$ and $X_\beta$.
4. Set $X_\alpha$ is called subset of set $X_\beta$ and is recorded as $X_\alpha \subseteq X_\beta$, if $\forall x (x[i_\alpha]x)_{Tr} \rightarrow (x[i_\beta]x)_{Tr}$. If at this $X_\beta \setminus X_\alpha \neq \phi$, we write $X_\alpha \subset X_\beta$.

It could be the situation, when sets $X_\alpha$ and $X_\beta$ do not have common points: $\overline{\exists} x : (x[i_\alpha]x)_{Tr} \wedge (x[i_\beta]x)_{Tr}$, i.e. overlap could be *null*; then we write $X_\alpha \cap X_\beta = \phi$.

To define plurality of sets $\{X_\alpha\}_{\alpha \in I}$ means to predetermine plurality of identifiers $\{i_\alpha\}_{\alpha \in I}$; here $I$ – is some indices, enabling us to perceive the differences between the identifiers, which we will always designate as by letter «i».

---

[1] Definition 4.1 differ from axiom ZF-V of *isolation* of sub-sections of ZF-theory particularly by analyzing identifier in the whole Universal of points, but not only in the "previously" isolated sets.

[2] *Axiom* ZF-VII. «There exists, at least, one set $X_\alpha$, having the following properties: a) $\phi \in X_\alpha$, b) if $x \in X_\alpha$, then $x \cup \{x\} \in X_\alpha$».
The first part of axiom has no sense for us, and the second one is tautology, since $x \cup \{x\} \equiv x$

[3] *Axiom* ZF-I. «The sets, consisting of the same elements, are equal».

**3.7. Comment.** But such *element* as "set $\mathsf{P}[X]$, the members of which are all kinds of subsets of set $X$" will not appear in RL-logic. Illusiveness of such structures of sets theory dramatically become apparent in identifiers medium: in the event of "all kinds of subsets", unification of any sub-plurality of identifiers coincides with some identifier from original plurality: $\{i_\alpha\}_{\alpha \in I} \equiv \mathsf{P}[\{i_\alpha\}_{\alpha \in I}]$. That is why for relations language Cantor theorem $\mathsf{m}[X] < \mathsf{m}[\mathsf{P}[X]]$ about potencies, where $\mathsf{m}[X]$ means potency of plurality $X$ in ZF in his understanding is unprovable in principle. Moreover, in the proof of this theorem (see [22], p. 29), when demonstrating the existence of inconsistent set, ill-defined statement $(x[r]\phi)$ is used ⌧

Hereinafter we will omit identifiers symbols in the instances, when there is be no special need for that.

The operations of overlapping, unification and taking the complement have the following, practically invisible properties, which are lumped in a form of theorem for convenience of references and with the purpose of completeness.

**3.8. Theorem.** Let $X,Y,Z$ be random sets.
1. $X \cup Y = Y \cup X$ и $X \cap Y = Y \cap X$.
2. $X \cup (Y \cup Z) = (X \cup Y) \cup Z$; $X \cap (Y \cap Z) = (X \cap Y) \cap Z$.
3. $X \cap (Y \cup Z) = (X \cap Y) \cup (X \cap Z)$; $X \cup (Y \cap Z) = (X \cup Y) \cap (X \cup Z)$.
4. $X \setminus (X \setminus Y) = X \cap Y$.

We also give without any proof de Morgan formulas, adduced below.

**3.9. Theorem.** $Y \setminus \cup \{X_\alpha : \alpha \in I\} = \cap \{Y \setminus X_\alpha : \alpha \in I\}$; $Y \setminus \cap \{X_\alpha : \alpha \in I\} = \cup \{Y \setminus X_\alpha : \alpha \in I\}$.

Asserting that set $Y$ exists in relation $r$ with set $X$, we will write $X[r]Y$. The plurality of those $y \in Y$, for which statement $(x[r]y)$ is true, is called *image* of relation $r$ in point $x \in X$, and is recorded as $r[x] = \{y \in Y : (x[r]y)_{Tr}\}$. The plurality of points out of $X$, for which $r[x]$ is not null, is called domain of *definition* of $r$, and $r[X]$ – *range of values* in $Y$.

Speaking, that set $Y$ is in relation $r$ with set $X$, we will imply that, if otherwise is not stated, that domain of $r$ definition is the whole set $X$.

And the last thing: relation $s$ is called inverse to $r$ and is designated as $r^{-1}$, if $\forall x,y \, (x[r]y)_{Tr} \square (y[s]x)_{Tr}$.

Let $r,s$ be some relations of set $Y$ with respect to set $X$. *Congruence* $r = s$ means $r[x] = s[x]$ for each point out of $X$; the same way *injection* $r \subset s$ is also defined, which implies $r[x] \subset s[x]$ in each point of the domain of definition (at least in one point $r[x]$ is proper subset of $s[x]$). In those cases, when injection $r \subset s$ *does not refuse* congruence, we will write $r \subseteq s$. *Unification* of two relations $r$ and $s$ is the name for relation $r \cup s$, the image of which in each point out of $X$ means just unification of images $(r \cup s)[x] = r[x] \cup s[x]$; the same way the *overlap* of relations $r \cap s$ is defined, as overlap of images of $r$ and $s$: $(r \cap s)[x] = r[x] \cap s[x]$.

Let $X[s]Y$ and $Y[r]Z$ be some relations between sets $X$, $Y$, $Z$. *Composition of* relations of $r$ and $s$ is the relation, designated as $r \circ s$, with image

$$(r \circ s)[x] = \cup \{r[y] : y \in s[x]\} \equiv \cup_{y \in s[x]} r[y],$$

the area of definition of which due to accepted agreements coincide with $X$; the composition is associative by definition. If $X = Y = Z$, then the composition is called commutative under condition $r \circ s = s \circ r$.

We totalize relations algebra, arising directly from the definitions, in the following form

**3.10. Theorem.** Let $r,s,t$ be relations, $X,Y$ – sets.
1. $(r^{-1})^{-1} = r$; $(r \circ s)^{-1} = s^{-1} \circ r^{-1}$.
2. $r \circ (s \circ t) = (r \circ s) \circ t$; $(r \circ s)[X] = r[s[X]]$.
3. $r[X \cup Y] = r[X] \cup r[Y]$; $r[X] \cap r[Y] \supseteq r[X \cap Y]$

We will distinguish between two types of relations. Relation $X[r]Y$ is called *algebraic* (mapping, translation invariant operator, or function – depending upon the context) and recorded in a form of r: $X \to Y$, if for each $x$ out of $X$ there exists unique element $y$ out of $Y$ and such, that $y = $ r$[x]$; and it is called multilingual, or *topologic* if not.

For two mappings there are specific names reserved. Mapping id$_X$: $X \to X$ is called *constant map* on $X$ (identical on $X$), if id$_X[x] = x$ for each $x \in X$; relation pr$_X$: $X \to Y$ is called *primitive* on $X$ (ill-defined on $X$), if pr$_X[x] = \phi$ for $x \in X$.

The notion of primitive relation is introduced in particular because in relations language unprovable statements explicitly appear. They arise also in those cases, when overlap of relations plurality, each of which has the same domain of definition, is primitive in each point out of the domain of definition: $\bigcap_\alpha$ r$_\alpha[x] = $ pr$[x]$.

Mapping f: $X \to Y$ is called *surjective* (mapping "on"), if f $\circ$ f$^{-1} = $ id$_Y$; *uinjective* (mapping "in"), if f$^{-1} \circ$ f $= $ id$_X$; *bijectie*, or one-to-one, if it is injective and surjective.

**3.11**. Statement. Let $X_0 \supset X_1 \supset X_2$ and f: $X_0 \to X_2$ be bijective mapping. Then, there exists bijection g: $X_0 \to X_1$.

*Proof.* In $X_0$ there exist subsets $Y_\alpha$ and such, that

$$Y_\alpha \supset \{X_0 \setminus X_1\} \cup \text{f}[Y_\alpha]. \tag{*}$$

The example could be the set itself $X_0$: $X_0 \supset \{X_0 \setminus X_1\} \cup \text{f}[X_0] = \{X_0 \setminus X_1\} \cup X_2$. Let us look at overlapping of *all* sets, complying with (*), i.e. $Y_0 = \bigcap Y_\alpha$. Then, antecedent (*) for this set reverts to congruence $Y_0 = \{X_0 \setminus X_1\} \cup \text{f}[Y_0]$, and in particular because there exist sets, for which injection (*) is realized in reverse direction; for example, for $X_2$ we have: $X_2 \subset \{X_0 \setminus X_1\} \cup \text{f}[X_2]$. Further, noticing, that by reason of 3.9 there is congruence $X_0 \setminus Y_0 = X_1 \setminus \text{f}[Y_0]$, we easily find required bijection g: $X_0 \to X_1$:

$$g[x] = \begin{cases} \text{id}[x], x \in X_0 \setminus Y_0 \\ \text{f}[x], x \in Y_0 \end{cases} \quad \boxed{}$$

This statement is instantly followed by famous theorem by Shroder-Bernstein, which plays important part not only in theory of sets, but also in theory of relations.

**3.12.** Theorem (Shroder-Bernstein). Let the following be in existence: bijective mapping p of set $X$ on some subset of set $Y$, and bijective mapping q of set $Y$ on some subset of set $X$. Then there exists bijection h: $X \to Y$.

*Proof.* Assuming $X_0 = X$, $X_1 = $ q$[Y]$, $X_2 = $ q $\circ$ p$[X]$, f $= $ q $\circ$ p, and using statement 3.11, we find the required bijection h: $X \to Y$ $\boxed{}$

For pre-images of mappings there exist certain correlations

**3.13**. Theorem. Let f: $X \to Y$ and $P, Q \subset Y$. Then
1. f$^{-1}[P \setminus Q] = $ f$^{-1}[P] \setminus $ f$^{-1}[Q]$.
2. f$^{-1}[P \cup Q] = $ f$^{-1}[P] \cup $ f$^{-1}[Q]$.
3. f$^{-1}[P \cap Q] = $ f$^{-1}[P] \cap $ f$^{-1}[Q]$.

*Plot* of mapping f: $X \to Y$ stands for set $Gr(\text{f})$, comprising pairs $(x, \text{f}[x])$. Any mapping is unequivocally characterized by its plot.

We will say, that mapping g, defined in some subset $P \subset X$, represents exposition of mapping f, if for all $x \in P$, g$[x] = $ f$[x]$ takes place. In this case, $Gr(\text{f}) \supset Gr(\text{g})$. In compliance with this, we will say, that f is *dilatation* of g.

Let $X$ be some set, and speaking about points $x, y, z,\ldots$ we will imply points out of $X$. Our attention will be concentrated on relation $X[\text{r}]X$.

**3.14**. Definition. Relation $X[\mathrm{r}]X$ is called
1. Reflexive, if $\forall x\,(x[\mathrm{r}]x)_{\mathrm{Tr}}$
2. Symmetric, if $\forall x, y\,(x[\mathrm{r}]y)_{\mathrm{Tr}} \Leftrightarrow (y[\mathrm{r}]x)_{\mathrm{Tr}}$
3. Transitive, if $\forall x, y, z\,(x[\mathrm{r}]y)_{\mathrm{Tr}} \wedge (y[\mathrm{r}]z)_{\mathrm{Tr}} \rightarrow (x[\mathrm{r}]z)_{\mathrm{Tr}}$
4. Equivalence relation, if it is reflexive, symmetric and transitive.

And such relations, starting from Section 4, will play a leading part, which rely on their well-known feature to break down the sets into non-overlapping *equivalence classes*.

**3.15.** Theorem. Let r be equivalence relation on $X$. Then, for any two points $x$ and $y$ sets $\mathrm{r}[x]$ and $\mathrm{r}[y]$ either do not coincide, or do not overlap.

*Proof.* Let us assume the contrary, that sets $\mathrm{r}[x]$ and $\mathrm{r}[y]$ do not coincide, but have one common point, let us say $z$. But then, due to transitiveness, any point out of $\mathrm{r}[x]$ also appertain to $\mathrm{r}[y]$, and consequently, $\mathrm{r}[x] \subset \mathrm{r}[y]$, which by virtue of symmetric property and reflexivity brings congruence. Hence $\mathrm{r}[x]$ and $\mathrm{r}[y]$ coincide, if they at least have one common point ⌛

If r is the relation of equivalence on $X$, we may consider new set $X^* = X/\mathrm{r}[\cdot]$, the elements of which are equivalence classes $x^* = \mathrm{r}[x]$ and which is called *factor set*, and the process of transfer to factor set itself is called *factorization* of $X$ with respect to r.

We should point out that, if r is equivalence relation on $X$, and for each other relation s, defined on $X$, $\mathrm{r} \subseteq \mathrm{s}$ is true, then we do not have any possibility to distinguish between point $x$ and its equivalence class $x^* = \mathrm{r}[x]$. That is why any object $x^*$ looks like point, until the relation is not found, which will "split" equivalence class into "different" points.

Not less substantial part in the future will be played by transitive-reflexive relations. Speaking about transitive-reflexive relations (in the event of need to point it out, we will write $\vec{\mathrm{s}}$), we will always imply, that it does not contain any equivalence relation, except for trivial one, or that respective set factorization is exercised with regard to this equivalence relation.

**3.16.** Definition. Set $X_{\mathrm{s}}$, with predetermined transitive-reflexive relation $\vec{\mathrm{s}}$ in it, is called
1. Partially s-ordered.
2. Normally s-ordered, if $\forall x, y\,\exists z : (\mathrm{s}[z] \subset \mathrm{s}[x])_{\mathrm{Tr}} \wedge (\mathrm{s}[z] \subset \mathrm{s}[y])_{\mathrm{Tr}}$.
3. Linear s-ordered, if $\forall x, y : (\mathrm{s}[x] \subset \mathrm{s}[y])_{\mathrm{Tr}} \vee (\mathrm{s}[y] \subset \mathrm{s}[x])_{\mathrm{Tr}}$.

In compliance with stated above, we will assert that "$y$ rank over $x$" and record $y > x$, if and only if $X_{\mathrm{s}} \setminus \mathrm{s}[y] \supset X_{\mathrm{s}} \setminus \mathrm{s}[x]$. Conformable to this we will speak about *maximal* element $\omega$ in $X_{\mathrm{s}}$, if $\mathrm{s}[\omega] = \{\omega\}$, and *minimal* element $\alpha$, if $\mathrm{s}[\alpha] = X_{\mathrm{s}}$. But if $\exists \kappa : \mathrm{s}[\kappa] \neq X_{\mathrm{s}}$, then $\nexists \alpha < \kappa$, and the plurality of such elements $\{\kappa\}$ will be called *root* elements (we have changed to some extent generally accepted terminology for more adequate to relations language).

We will call two random elements $x, y \in X_{\mathrm{s}}$ out of partially s-ordered set as *comparable*, if $(\mathrm{s}[x] \subset \mathrm{s}[y])_{\mathrm{Tr}} \vee (\mathrm{s}[y] \subset \mathrm{s}[x])_{\mathrm{Tr}}$; and respectively non-comparable, if $(\mathrm{s}[x] \not\subset \mathrm{s}[y])_{\mathrm{Tr}} \wedge (\mathrm{s}[y] \not\subset \mathrm{s}[x])_{\mathrm{Tr}}$.

Maximal elements, as well as root ones, could be plenty of, and all of them are non-comparable.

**3.17**. Statement. Minimal element, if it exists, is always singular.

*Proof.* Really, let $\alpha$ and $\alpha'$ be two minimal elements. Hence $\mathrm{s}[\alpha] = X_{\mathrm{s}}$ and $\mathrm{s}[\alpha'] = X_{\mathrm{s}}$. But then $\alpha' \in \mathrm{s}[\alpha]$ and $\alpha \in \mathrm{s}[\alpha']$, and that is why transitive-reflexive relation s is symmetric in certain subset $X_{\mathrm{s}}$, and consequently, incorporates nontrivial equivalence relation, which is contrary to definition of transitive-reflexive relation ⌛

And here comes the last out of anticipatory definitions and theorems.

**3.18**. Statement. For any plurality of mutually overlapping sets $\{X_\alpha\}_{\alpha \in I}$ there exists mapping $\mathrm{f} : I \to \bigcup X_\alpha$, which is called *choice function*, and it is such, that $\forall \alpha\,\mathrm{f}[\alpha] \in X_\alpha$.

*Proof.* Antecedent $x \in X_\alpha$ matches $(x[i_\alpha]x)_{\text{Tr}}$, and therefore $i_\alpha = \bigcup_{x \in X_\alpha} i_{\alpha,x}$, where $i_{\alpha,x}$ is defined by the antecedent, that $\exists! x \in X_\alpha : (x[i_{\alpha,x}]x)_{\text{Tr}}$. Now for each $\alpha \in I$ we assign some $i_{\alpha,x_\alpha}$ – it exists, though, maybe, it is not singular, and the only remaining thing is to assume

$$\forall \alpha \; f[\alpha] = \{x : (x[i_{\alpha,x_\alpha}]x)_{\text{Tr}}\} \equiv x_\alpha \in X_\alpha \; \boxtimes$$

In mathematical literature statement 3.18 is called *axiom of choice*. For relations language the issue of choice function is the direct consequence of the point identifier *definition*. ; but the *choice* has to be done. That is why we will still call 3.18 as the axiom of choice.

**3.19**. Comment. Axiom of choice very often, though tacitly, was used in mathematics. But, after it was formulated by Zermelo, the mathematicians' community split down the middle. Some started accepting this axiom without any doubt, the other started to reject the same way. And the reason was that on its basis Zermelo proved, that arbitrary set $X$ could be easily *put in order* – could be put in order in such a way, that each subset will have minimal element. In other words, it is considered proved in the theory of sets, that choice axiom matches, as we will say at present, Cantor *axiom* [1]:

**1**. Axiom (Cantor). For random set $X$ *there exists* transitive-reflexive relation s on $X$ and such, that $\forall Y \subset X \; \exists y \in Y : s[y] \supseteq Y$.

Further, complete order is easy to construct: from 3.19.1 it flows out: $\exists x_1 : s[x_1] = X$, and by reason of definition of transitive-reflexive relation, this element is singular and minimal. Let us consider now subset $X \setminus x_1 \subset X$. Again, due to 3.19.1 and 3.17 we easily find:
$\exists! x_2 : s[x_2] = X \setminus x_1$. And also $x_1 < x_2$. Continuing by induction: $\exists! x_\mu : s[x_\mu] = X \setminus \{\bigcup_{\eta < \mu} x_\eta\}$, ...and stating, where necessary, the existence of "actually infinite number", we get, that any set could be easily put in order.

But the key element of equivalency proof of choice axiom and axiom 4.19.1 (see [23], p. 28; [24], p. 78) rely on the existence of particularly such choice function, for which $f[\phi] = x$. And this is ill-defined statement. Ill-defined statement of such kind is tacitly present also in the other methods of Zemerlo theorem proof (see [25], p. 82-83).

That is why we formulated Zemerlo *theorem* in the form of Cantor *axiom* – since it is particularly this thesis, which Cantor failed to prove, (and which he considered to be "cogitation law"), which ran through all the papers dedicated to theory of sets before Zemerlo's published paper.

Putting it the other way, it was not the axiom itself, which agitated mathematicians community (see, for instance, [16, 26] and the literature cited there), but "materialization" of null-set, tacitly concealed in the proof procedure. Moreover, axiom 3.19.1, to our mind, can not match choice axiom in concept – it incorporates two coherent existential quantifiers, precisely:

**2**. Axiom (Cantor). $\forall X \; \exists \vec{s} : [\forall Y \subset X \; \exists y \in Y : \vec{s}[y] \supseteq Y]$,

and statement 3.18**,** being analyzed as axiom, has only one.

But Cantor axiom can not also be invalidated by simply referring to impotence of our mind, which just is not able to find the respective transitive-reflexive relation.

It is important to stress, that Cantor axiom makes it possible to prove the "principle of transfinite induction", which plays such substantial part in the theory of sets, that one can not simply construct the theory of sets in the pure meaning of the word.

**3**. Theorem. Let $X_s$ be perfectly s-ordered set, each element of which is assigned to statement $p[x]$, from truthfulness of which for all $x < y$ it arises, that $p[y]_{\text{Tr}}$. Then, it flows out, that $p[x]_{\text{Tr}} \; \forall x \in X_s$.

*Proof.* Let us assume the opposite, that set $Y = \{y \in X_s : p[y]_{\text{Li}}\}$ is not null. But it has minimal element, which we designate as $y_\alpha$. Then for all $x < y_\alpha \; p[x]_{\text{Tr}}$; and this brings $p[y_\alpha]_{\text{Tr}}$. Consequently $Y = \phi \; \boxtimes$

---

[1] «The notion of *well ordered set* proves to be fundamental for the whole teaching of manifolds. In one of the coming papers I will return to what seems to me to be principal, open to many hazards and especially distinctive by its general significance – cogitation law, in compliance with which each *strictly defined* set could be shaped as *well ordered* set». Quoted from [22], p. 69. Italics are Cantor's.

That is why in general case this theorem is also unprovable. But in the set of whole numbers and their Cartesian products in an explicit form – see 3.20.4, it is possible to indicate transitive-reflexive relation and of such kind, that each subset will have minimal element ⌛

Hereinafter we will consider several examples of differently ordered sets.

**3.20.** Examples.

1. In set of congruums **N** there exists natural order, connected with definition of whole numbers – see below. We will designate transitive-reflexive relation of natural order as $e$. In point $n \in \mathbf{N}$

$$e[n] = \{m \in \mathbf{N}: m \geq n\}.$$

Relying upon natural order, connected with definition of congruums and negative integers, let us predetermine the following transitive-reflexive relation $s$ in **N** with image

$$s[n] = \{m \in \mathbf{N}: (-1)^n m \geq (-1)^m n\}.$$

It is easy to deduce, that $s[0] = \mathbf{N}$: and $s[1] = \mathbf{1}$, and the rest of the numbers, at first even in natural order, and later odd in reverse order, are located between these minimal and maximal elements: 0, 2, 4, 6, … …7, 5, 3, 1.

Set $\mathbf{N}_e$ is quite ordered, and $\mathbf{N}_s$ is not: subset of odd numbers does not have minimal element. We can try to introduce $s$-minimal of all odd numbers, but $s$-maximal of even numbers is actually infinite number $\omega$, which isolates even and odd numbers in $\mathbf{N}_s$. But in order to prove, that $(-1)^\omega = 1$, we will have to use the principle of transfinite induction, which rests just on existence of such order in set, that each subset has minimal element. And this means circulus vitiosus.

2. Let $\Sigma_s = \{x_1, x_2, x_3, x_4\}$ and $s[x_1] = \{x_1, x_3, x_4\}$, $s[x_2] = \{x_2, x_3, x_4\}$, $s[x_4] = \{x_4\}$.

We see partially $s$-ordered set, without minimal element, with two maximal and two root elements. Pairs of elements $x_1, x_2$ и $x_3, x_4$ are incommensurable.

3. Now $E = [0,1]$ is interval of real numbers, $Y = E \times E$ and $s[x_1, x_2] = \{(y_1, y_2): x_1 \leq y_1, x_2 \leq y_2\}$. Transitive-reflexive relation $s$ predetermines normal order in $Y$: for each two different points $x = (x_1, x_2)$ and $y = (y_1, y_2)$ there exists $z = (z_1, z_2)$ and of such type, that $s[z] \subset s[x]$ and $s[z] \subset s[y]$ – Fig.1. Later, $s[0,0] = Y$, and $s[1,1] = (1,1)$.

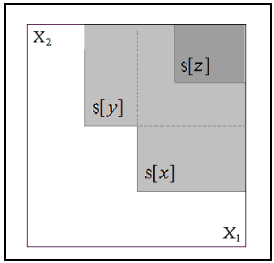

**Fig.1**. Schematic representation of image $s[\cdot]$ in different points – darkened zones, including boundaries.

4. Let us obvert square $Y = E \times E$ into linearly ordered set [25]. We assume, that $r = s_1 \bigcup s_2$, where $s_1[x_1, x_2] = \{(y_1, y_2): x_1 < y_1, y_2 \in E\}$, $s_2[x_1, x_2] = \{(y_1, y_2): x_1 = y_1, x_2 \leq y_2\}$.

For two different points $x = (x_1, x_2)$ and $y = (y_1, y_2)$ it is true, that either $r[x] \subset r[y]$, or $r[y] \subset r[x]$ – Fig.2. And also $s[0,0] = Y$ and $s[1,1] = (1,1)$.

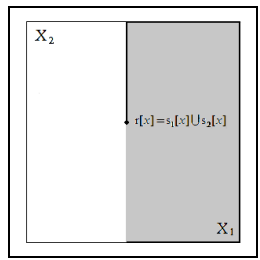

**Fig.2**. Schematic representation of image $r[x]$, which is darkened zone, including also heavy line with the point.

If **N** is selected as co-ordinate sets, then analyzed relation linearly puts in order $\mathbf{N} \times \mathbf{N}$, and in each of its subsets minimal element will be found ⌛

With respect to whole numbers (see below about the numbers), the sets in the theory of relations are subdivided into three categories: *finite, denumerable,* or *non-denumerable*. But, first we should give one convenient designation: $N|_n = \{1, 2, ..., n\}$.

**3.21**. Definition.

1. Plurality $X$ is called finite (infinite), if whole number $n$ is found (not found), and such number, that bijection $f : X \to N|_n$ exists.

2. Infinite plurality $X$ is called denumerable (non-denumerable), if bijection $f : X \to N$ exists (does not exist).

Thereinafter we will bring forward a range of theorems, which are not only instructive, but which will be used in the future. Those theorems are associated with the notion of denumerability, along with respective changes, complying with the relations language.

**3.22.** Theorem. Denumerable aggregate of non-crosscutting denumerable sets is denumerable.

*Proof.* Let $\{X_n\}$ be denumerable plurality of non-overlapping denumerable sets. We assume, that $X_n = \bigcup_k x_{n,k}$. Then, their aggregate can, for example, be recorded in the following form

$$\bigcup X_n = \{x_{1,1}, x_{2,1}, x_{1,2}, x_{3,1}, x_{2,2}, x_{1,3}, ...\}$$

($n + k = m$; $m = 2, 3, 4, ...$), or in the other one, which is analogous. Now one-to-one accordance with N is compelling ∎

**3.23.** Theorem. In each non-denumerable set $X$ there exists denumerable subset $X_N$.

*Proof.* Let us assume that $X_1 = X$. By virtue of 3.18 there exists choice function $\mathbf{f}_1$, and such function, that $\mathbf{f}_1[1] = x_1 \in X_1$. Let us consider set $X_2 = X_1 \setminus \{x_1\}$. Further, again by virtue of 3.18, there exists choice function $\mathbf{f}_2$, and such function, that $\mathbf{f}_2[2] = x_2 \in X_2$. Continuing in inductive manner, we will get injective mapping $\hat{g}: N \to X$, which is being given by rule: $\hat{g}[n] = \mathbf{f}_n[n]$. Image $\hat{g}[N] \subset X$ is a denumerable subset ∎

**3.24.** Theorem. Let $X$ be non-denumerable, and $Y_N$ – denumerable sets. There exists bijection $f : X \bigcup Y_N \to X$.

*Proof.* Because of 3.23 let us isolate denumerable subset $X_N$ out of $X$. Then

$$X \bigcup Y_N = (X \setminus X_N) \bigcup X_N \bigcup Y_N = (X \setminus X_N) \bigcup (X_N \bigcup Y_N)$$

Because of 3.22 there exists bijection $g : X_N \bigcup Y_N \to X_N$ and identical self-mapping $X \setminus X_N$. That is why desirable bijection is the following: $f[x] = \begin{cases} \text{id}[x], x \in X \setminus X_N \\ g[x], x \in X_N \bigcup Y_N \end{cases}$ ∎

**3.25.** Theorem. Let $X$ be non-denumerable, and $Y_N$j – denumerable sets. There exists bijection $f : X \to X \setminus Y_N$.

*Proof.* Set $X \setminus Y_N$ is non-denumerable, since if not, because of 3.22, set $(X \setminus Y_N) \bigcup Y_N = X$ would have been denumerable. But from 3.23 it flows out, that there exists one-to-one mapping $f : (X \setminus Y_N) \bigcup Y_N \to X \setminus Y_N$, which precisely represents desirable ∎

**Products of sets**

Now we would like to say some words about very important construction of the theory of sets – sets product. This construction has already been tacitly used by us for the examples 3.20.3,4, and we will investigate it in details in future.

Let $\{X_\alpha\}_{\alpha \in I}$ be some kind of sets plurality, where $I$ is some plurality of indices. In the case of denumerable set $I$ *product of sets* is defined univalently

$$Y = X_{\alpha_1} \times X_{\alpha_2} \times X_{\alpha_3} .. \equiv \prod_{\alpha \in I} X_\alpha,$$

as plurality of, let us say, "words" $(x_1, x_2, x_3, ..)$:

$$\prod_{\alpha \in I} X_\alpha = \{(x_1, x_2, x_3, ..) : x_1 \in X_{\alpha 1}, x_2 \in X_{\alpha 2}, x_3 \in X_{\alpha 3} ..\},$$

where conditions of "letters" $x_\alpha$ belonging to *co-ordinates sets* $X_\alpha$ and their order are *enumerated*. Generalization of co-ordinate sets for *non-denumerable* set is possible, but only in the case, if indices set $I$ is *perfectly well ordered*. It is only then, when we can say, that co-ordinate set $X_\beta$ *is a successor of* $X_\alpha$, or quite the reverse, and put down "word" $(x_1 x_2 x_3 ... x_\alpha x_\beta ...)$, which complies with the point of product.

At each $\beta \in I$ mapping of projecting is defined $p_\beta : \prod_{\alpha \in I} X_\alpha \to X_\beta$ for $\beta$-multiplier of product according to the principle: $p_\beta[x] = x_\beta$ for each $x \in \prod_{\alpha \in I} X_\alpha$.

### About numbers and their notation

The notion of whole number "one", for which graphic symbol «1» is adopted, presents primary notion. It is pure "form" without reference to substance, and consequently, it also presents *equivalence relation*, explicitly expressed by this notion. Claims of the theory of sets to definition of the whole number as the category of "equivalent sets" are groundless, since in the notion of equivalent sets the notion of bijection – **one**-to-one correspondence is already imbedded, and consequently the notion of "one".

The whole numbers are defined inductively. For example, record type «2 = 1+1» is *definition* of graphic symbol «2», and etc. The whole numbers, along with the principle of mathematical induction, into which the principle of transfinitary induction degenerates in normally ordered set of the whole numbers, precede [4] any science. Poincare used to write: *"...it is difficult to say a phrase, without introducing into it word «number», or word «several», or, after all, some word in plural number"*. And further, not without poignant sarcasm: *"Burali-Forti defines number 1 in the following way:* $1 = iT'\{Ko \cap (U,h)\varepsilon(U\varepsilon, Un)\}$. *This definition is extremely suited to provide a glimpse of number 1 to those, who have never heard anything about it!"*. And then he points out circulus vitiosus both in this "definition", and in all subsequent attempts to define the whole numbers proceeding from "arche-principles".

Subtle representation of definition of the whole numbers can be found in the first chapters [4,13]; we will not repeat ourselves. But we would like to consider the issue of "presentation of real numbers" and proofs of some theorems associated with it.

Let $E = [0,1]$ be unit segment of real numbers, $X_n = \{0, 1, ..., n-1\}$ – plurality, comprised of whole numbers $0, 1, ..., n-1$, and $D_n = \prod X_n$ – denumerable Cartesian product of such sets, the element of which "word" $(x_1, x_2, ...)$, where any "letter" $x_i$ is equal to certain number out of $0, 1, ..., n-1$. Index of the letter is called place of the letter.

Mapping $f : D_n \to E$, which puts in correspondence to each word $(x_1, x_2, ...)$ out of $D_n$ number $x$ out of $E$ in compliance with convention $x = \sum x_i / n^i$, is surjective – common number can correspond to two different words. In order to make the correspondence bijective, the agreement is reached to reject the words, having either $0$, or $n-1$ in units period. Most often the rejected words are those, which have 0 in their period. That is why we will agree to consider words $(x_1, x_2, ..., x_k(0))$ "dead", whereas the rest of the words – "alive", also including word $(0(0))$, which has unique representation. Let us designate the set of "dead" words as $M_n = \bigcup (x_1, x_2, ..., x_k(0))$. Now the abovementioned mapping $f : D_n \setminus M_n \to E$ is bijective.

It is easy to see, that the number of "dead" words is denumerable: the plurality of words with length up to $k$-steenth place is finite, and $M_n$ presents their countable union. And further on – 3.22.

Underneath we will consider five theorems, this or that way connected with Cantor name, and which once laid the foundation for active investigations in the theory of sets.

**3.26.** Theorem. Set $E$ is non-denumerable.
*Proof* 1. In proving non-denumerability of real numbers set $E$, the opposite is assumed, that set of alive words, univalently cohering $E$, is denumerable. And in $n$-fold notation, graph of alive words is constructed

$$x_1 = (x_{1,1}, x_{1,2}, ..., x_{1,k}, ..)$$
$$x_2 = (x_{2,1}, x_{2,2}, ..., x_{2,k}, ..)$$
$$...............................$$
$$x_n = (x_{n,1}, x_{n,2}, ..., x_{n,k}, ..)$$
$$...............................$$

Later, by "Cantor diagonal method", word $y = (y_1, y_2, ..., y_k, ..)$ is constructed, where "$k$-*steenth letter $y_k$ of word $y$ is selected different from letter* 0 *and letter* $x_{k,k}$ *of word* $x_k$". As a result of it, we get alive word, which do not belong to the initial denumerable set, since it differs from any other word from the graph by some kind of letter. It is contradiction, which proves, that the set of alive words can not be denumerable. The key role here is played by "magic phrase", italicized above, which is three-pace function $y_k = p[0, k, x_{k,k}]$, and not letters selection function f, which, by reason of its definition, can not be determined by merely number of the word or column in the graph. If, for example, we assume $y_k = (n-1) - f[k]$, where $f[k] = x_{k,k}$, then by virtue of *graph equivocation*, disregarding what kind of function the choice function could be – not necessarily $f[k] = x_{k,k}$, it may turn out, that, starting from some $k > m$, it will be precisely $(n-1) - f[k] = 0$. And the obtained word will belong to the set of dead words, just because of that it does not appertain to source graph. Hence, there is no contradiction.

*Proof* 2. The other one, which also ascends to Cantor, among the existing proofs of non-denumerability of unit segment points set, was mastered by us in Section 1. It is barely fallible.

*Proof* 3. Later, in theorem 3.27 we will see, that there exists bijection $g:D_2 \to E$, though not so trivial, as for $f:D_2 \setminus M_2 \to E$. And we just have to convince ourselves, that set $D_2$ is non-denumerable. But again, assuming the opposite, that $D_2$ is denumerable set, we construct graph of words out of $D_2$, and by "Cantor diagonal method" the word $y = (y_1, y_2, ..., y_k, ..)$ is constructed, where $y_k = 1 - x_{k,k}$, which does not appertain to source graph of words. We pay attention that, at stated construction of word, choice function is strictly defined: $\forall k \; f[k] = x_{k,k}$ ▆

**3.27. Theorem.** There exists bijection $g:D_2 \to E$

*Proof* 1. Relying upon denumerability of set $M_n$ and theorem 3.25, existence of bijection $g:D_n \to E$ is compelling.

*Proof* 2. Let $k/n$ be any irreducible fraction out of whole numbers $k \geq 2$ and $k/n < 1$. Let us consider subset $K_{k/n}$ out of $E$, comprising numbers category $x = \sum k x_i / n^i$, where word $(x_1, x_2, x_3, ...)$ is taken out of $D_2$. Taking into account, that mapping $f:D_2 \setminus M_2 \to E$ is bijective, we get the following diagram of associations

$$\begin{array}{ccc} D_2 & \xrightarrow{h} & K_{k/n} \\ \cup & & \cap \\ D_2 \setminus M_2 & \xleftarrow{f^{-1}} & E \end{array}$$

As a result of it: set $E$ is bijectively mapped on subset $D_2 \setminus M_2$ of set $D_2$, and set $D_2$ is bijectively mapped on subset $K_{k/n}$ of set $E$. And Shroder-Bernstein theorem predicates existence of bijection $g:D_2 \to E$ ▆

From proof 2, as direct consequence we get

**3.28. Theorem.** There exists bijection $h \circ g^{-1}: E \to K_{k/n}$.

**3.29. Comment.** Plurality $K_{2/3}$ matches Cantor perfect set. We have to remind, which way it is comes. Unit segment is divided into three equal parts, and the middle part is cleared out. The remaining segments are again divided into three parts, and in each of them the middle interval is cleared out, and etc. The thing, which remains after such calculation procedure is called Cantor perfect set. At this, the sum of lengths of the cleared out intervals is equal to one. It was the first example of geometrically intuitively constructed *non-demunerable*, equivalent to unit segment of pretty "thin" set ▆

**3.30. Theorem.** There exists bijection $f: E^m \to E$.

*Proof.* First of all, let us consider the proof of existence of bijection between unit segment $E$ and its square $E^2 = E \times E$.

The proof rests on Cantor observation, that for each two words $(x_1, x_2, ...)$ and $(y_1, y_2, ...)$ compound word $(x_1, y_1, x_2, y_2, ...)$ also turns to be one-to-one. If the original words were alive, then the compound word would for sure be alive. However, the opposite assertion will be l. For example, any alive word type $(x_1, x_2, ..., x_k (01))$, or $(x_1, x_2, ..., x_k (0001))$ and etc., can not be represented by two alive words, and hence, as long ago as Dedekind noticed, they do not have their representation in $E \times E$. The above-stated does not depend on the numbers notation, which we choose.

However the required bijection exists. Noticing one-to-one mapping of the interval in square subset, for example $g: E \to E \times (0)$, and taking into consideration the above bijection of the square in the interval subset, the only remaining thing is to use Shroder-Bernstein theorem.

The result is plainly generalized by any *finite* product $E^n$. But even for denumerable product it is impossible to prove by the above procedure the existence of one-to-one mapping $E^N$ in subset $E$, since we will not be able to record "word", matching any point $E$, because the whole "length" of word $E$ will be occupied by the first three letters of the words out of $E^N$ ⌛

Let $F_N$ be a designation for the set of all functions, mapping the set of positive numbers $N$ in two-point set $X_2 = \{0, 1\}$.

**3.31.** Theorem. There exists bijection $f: E \to F_N$.

*Proof.* For each $g \in F_N$ the consecutive order of its value on $N$ represents consecutive order of zeros and ones, and consequently each function in one-to-one manner corresponds to some point out of $D_2$. But, as we have seen in 3.27, there exists bijection $g: D_2 \to E$, and that is why bijection $f: E \to F_N$ does exist ⌛

Let $F_E$ designate the plurality of all functions, mapping $E$ in $X_2 = \{0, 1\}$.

**3.32.** Theorem. Bijection $f: E \to F_E$ does not exists.

*Proof.* Let us assume the opposite – that such bijection $f: E \to F_E$ exists. It is choice function $f[x] = g_x$ and by assumption it is of the kind, that $F_E = \bigcup_{x \in E} g_x$. However, function $h[x] = 1 - g_x[x]$ does not appertain to this set, since it differs from any of the functions $g_x$ by its value precisely in point $x$. Consequently, there is no bijection.

But there exists bijection $E$ in subset $F_E$: for each point $x$ function $f_x[y]$ is put in correspondence, and the function takes, for instance, value 1 in point $y = x$, and 0 in other points ⌛

**Comparative survey**

Hereinafter we will carry out brief comparative survey of the basic notions of ZF-theory of sets, see, for example [16], which, for the sake of picturesque, we will call "the world of Georg Cantor", and RL-theory of sets, which we will, respectively call "the world of Henri Poincare".

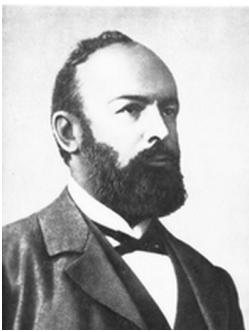

**Photo 1**.
Georg Cantor
(1845-1918)

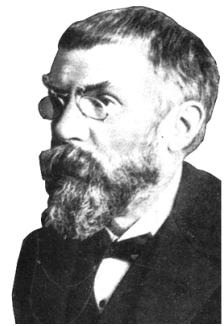

**Photo2**.
Henri Poincare
(1854 – 1912)

In the world of Poincare *Universal of points* is described – immense object as well as the relations, in which the points exist with respect to each other. Each of the points is an undefinable object of our consciousness. It is a mental image of some relation of equivalency. And the "point" will remain the "point", until the relation is found, which will "split" it into all kinds of "points". That is why the notion "potency of a set" for relations language looses its sense – it becomes "nondimensionalized" with respect to "language luxuriance", used in the theory, and this emphasizes relativism of notions of sets theory.

In Universal new points can also appear, when new relations arise, which bind this point with the other points. Before that it would be beyond the "visible horizon" of the Universal. That is why the horizon of points Universal moves far and far away, while we cognize it deeper and deeper.

Sets in the world of Poincare are nothing more than "distinguished" – "categorized", or "filtered" points of Universal; they are not singled out by set name. If set is infinite, then with respect to each "produced" point we can always say, whether it belongs to the set under consideration, or not. But the infinite set itself in the world of Poincare is impossible to "produce" – it is always just "potentially" infinite sets.

In the world of Cantor there exists *Universal of sets*. The horizon of this Universal, even if it is fixed, it is fixed only under the antecedent that the Universal itself is other than a set. And the set is not just "plurality" of points, which can be conceived as single entity, but precisely as entity, which will obligatory be at the same time a "member" of some other set. Moreover, infinite sets in the world of Cantor always represent "actually" infinite sets.

It is a kind of visionary world: "ex nihilo" – out of "null-set" $\phi$, "something" – $\{\phi, \{\phi\}, \{\{\phi\}\}, ...\}$ is easily obtained; "cogitation" is sufficient to represent the set as a point of the other set.

As a result, there appear such objects, as "set $\mathsf{P}[X]$, the elements of which are all kinds of subsets of set $X$". And there are no limits for such objects, as $\mathsf{P}^n[X]$ and etc. It makes it possible within the course of proving some theorems to view different, even not overlapping subsets, either in the quality of different points (and then, for example, mapping $f:\mathsf{P}[X] \to Y$ could be found), or again in the quality of subsets. Maybe, it is just this duality of the notion of set-object, where readiness lies in, which facilitates overcoming of difficulties, associated with proofs of some theorems.

And for relation language null-set does not exist: it is invariably that $x \equiv \{x\}$, and the set appears as a point only when it is brought to life by some relation of equivalency; now, there exist sets, for which $X = \mathsf{P}[X]$, but always $\mathsf{P}^2 = \mathsf{P}$; in an explicit form unprovable statements appear, etc.

If in the world of Poincare the primary interest particularly lies in investigation of relations among the sets, then in the world of Cantor the relations play though very important, but all the same, auxiliary role. Among the possible inter-set relations the most important is one-to-one correspondence – bijection. In Universal of sets bijection is relation of equivalence. That is why the equivalent class is *objectivized* and named *potency* of a set.

But to our point of view, the whole world of Cantor rests on one fundamental premise – on Cantor axiom 3.19.1. This axiom, never decisively asserted by Cantor, is not "the law of ideation", but namely *axiom of existence* in random set of certain transitive-reflexive relation, which together with Shoeder-Bernstein theorem guaranteed for any two sets, that either one of them is bijectively mapped in some subset of the other, or they are equivalent. That is why subsets potencies can be compared as usual numbers. Afterwards confidence in this was given by Zemerlo theorem, which narrowed down Cantor axiom to axiom of choice, the axiom, which before that moment seemed so trivial, that it was not clearly stated. Then in reality there appeared such thing, which now is called theory of sets.

Except that, neither Cantor theorem, guaranteeing the existence of *potencies scale*, nor Zemerlo theorem, guaranteeing the possibility to absolutely put in order any set, can not be considered proved for relations language, since they use indefinite statements. Moreover, Cantor axiom, to our mind, in principle can not be narrowed down to choice axiom: it incorporates two linked existential quantifiers, and choice axiom – only one.

On the other hand, it is possible just to accept Cantor axiom – it can not be refuted; especially because there exist examples for which it is fulfilled. In regard to potencies scale, its "existence" is easy to be grounded on the examples of existence of denumerable, non-denumerable and non-equivalent sets, for instance Cantor theorem 3.23. And since these "potencies" can be compared, then, maybe, these are the two circumstances, which once upon a time compelled Hilbert to declare, that: *"Nothing can exostracize us from paradise, created for us by Cantor"*. And it is really so.

In the world of Poincare, due to indeterminateness of Universal boundaries, there is no such a notion as "sets equivalence class" – it is simply not defined, and consequently there is no notion of "potency". On the other hand, the notion of equivalence of two sets makes it possible, in concept, to speak, following Neyman, about "representer" of this "class". But the relations, in which these "representers" exist is determined particularly by Cantor axiom and transfinite induction theorem, flowing out of it.

But it is just Cantor axiom that we do not accept. To our mind, this axiom not only narrows the range of vision (as we already pointed out many times, it contains two linked existential quantifiers, and statement 3.18, being considered as choice axiom, has only one), but, and it is the most important thing, via the null-set, it opens the door into the world of statements, which are undefined and consequently unprovable by relations language.

And the *existence* of choice function in the world of Poincare is accepted not as "choice axiom", but it is the direct consequence of point identifier definition, though at any specific construction the *choice*, as such, has to be made.

As in the world of Cantor, the notion of bijection is one of the most important notions in the world of Poincare. But now the issue of existence of bijection between the sets, or, speaking the language of theory of sets, definition, whether these sets are located in the same "class", or in which relations "class represeners" exist, each time requires discrete exploration. We set enough examples, moreover classical ones, which way bijection absence or existence is determined, without resorting to the notion of cardinalities. For sure, it gives hard time, makes the life more prosy, but in return, saves us from all and even potential paradoxes of sets theory. It is just because in the world of Poincare there are no indefinite statements, the relations are specified, and …

And the only thing remaining is just to "calculate". Maybe, all these things are a bit sobersided, but it is the reality of the world of Henri Poincare. And though as the result of calculations we can come to an unprovable statement, it should not embarrass us. We are not gods, at long last, as it seemed to someone at the beginning of the XX century.

## Section 4. Binary relations filters

What is usually meant by the words "mathematical structure"? As Bourbaki used to write [27], "*The common feature of different notions, united by this generic term, is their applicability to the set of elements, the nature of which is not found. In order to identify the structure, one or several relations, in which its elements exist, are preset; then it is stated, that the given relation or relations satisfy some conditions(which are specified and which are axioms of the structure under consideration)*".

But in this interpretation of *structure* word-combination "some conditions" provides too wide latitude for whims – any objects-points satisfy *some* conditions. But everybody understands that a bourock is not a building – a synonym of the word "structure". The relations, in which the sets exist, should contain certain common ab indication that we really deal precisely with a building, but not a bourock.

So, any mathematical structure represents subset $X$ of elements of some set $Y$, the subset isolated by some characteristics.

Actually, this way we defined the notion of "set" $X_\alpha$, where in the quality of $Y$ we had Universal of points $F$, and a the indicator of point $x$ appurtenance to $X_\alpha$ was identifier $i_\alpha$ of set $X_\alpha$, which reproduced only the property of *reflexivity* of relation $i_\alpha$ on Universal of points $F$.

Adding to the property of relation to be reflexive the property to be also *transitive*, we come to "ordered structures" of four types: *partially*, *normally*, *linearly*, and *well-ordered* sets.

In regard to the structures *univalently* defined by any of its elements and brought to existence by some plurality $I$ of relations $\{r_\alpha\}_{\alpha \in I}$, whether it be a group, topological space and etc, it is necessary to require from the relations, which induce the structure, to involve combined properties of reflexivity, transitivity and symmetry. It is only then, that a structure is induced by some relation of equivalency, let us say, $J_r$ in $F$, expressing some kind of typical features. But relation of equivalency $J_r$ is not known to us. But we should note, that if two points out of $F$ exist in certain relation $r_\alpha$, they should also exist in relation $J_r$; and the opposite, if a pair of points exist in relation $J_r$, then they must exist in certain relation $r_\alpha$ out of plurality. By this required $J_r$ emerges as the smallest upper edge of relations plurality $\{r_\alpha\}_{\alpha \in I}$, which induces mathematical structure, i.e. $J_r = \bigcup_{\alpha \in I} r_\alpha$; accordingly, the largest lower edge will be designated as $j_r = \bigcap_{\alpha \in I} r_\alpha$.

Both the lower and the upper edges can belong to the plurality. Then we may speak about minimal and, respectively, about maximal elements of the plurality. And in order to obtain the description of all the mathematical structures, univalently defined by any of its elements, the only thing remaining is to find out the conditions, which should be imposed on elements of relations plurality, for $J_r$ to be the relation of equivalency.

Hereinafter, in three definitions we lumped different conditions for the elements of relations plurality, which ensure the property of the upper edge to be reflexive, symmetric and transitive. The conditions are put down in such a way, that satisfaction of any of the axioms results in satisfaction of the axioms, which are predecessors of it.

**4.1. Definition.** Plurality of Relations $\{r_\alpha\}_{\alpha \in I}$ in $F$ is called
1. Weakly reflexive, if $\forall x \exists r_\alpha : (x[r_\alpha]x)_{Tr}$.
2. Normally reflexive, if $\exists r_o : \forall x \, (x[r_o]x)_{Tr}$.
3. Strongly reflexive, if $\forall r_\alpha \forall x \, (x[r_\alpha]x)_{Tr}$.

**4.2. Definition.** Plurality of relations $\{r_\alpha\}_{\alpha \in I}$ in $F$ is called
1. Weakly symmetric, if $\forall r_\alpha \forall x, y \exists r_\beta : (x[r_\alpha]y)_{Tr} \to (y[r_\beta]x)_{Tr}$.
2. Normally symmetric, if $\forall r_\alpha \exists r_\beta : r_\beta \supset r_\alpha^{-1}$.
3. Uniformly symmetric, if at normal symmetric properties $\forall r_\beta \exists r_\alpha : r_\beta \supset r_\alpha^{-1}$.
4. Strongly symmetric, if $\forall r_\alpha : r_\alpha = r_\alpha^{-1}$.

**4.3. Definition.** Plurality of relations $\{r_\alpha\}_{\alpha \in I}$ in $F$ is called
1. Weakly transitive, if $\forall r_\alpha, r_\beta \forall x, y, z \exists r_\gamma : (x[r_\alpha]y)_{Tr} \wedge (y[r_\beta]z)_{Tr} \to (x[r_\gamma]z)_{Tr}$.
2. Normally transitive, if $\forall r_\alpha, r_\beta \exists r_\gamma : r_\gamma \supset r_\alpha \circ r_\beta$.
3. Locally uniformly transitive, if at normal transitivity there is also $\forall x \forall r_\gamma \exists r_\alpha, r_\beta : r_\gamma[x] \supset r_\alpha \circ r_\beta[x]$.
4. Uniformly transitive, if at normal transitivity $\forall r_\gamma \exists r_\alpha, r_\beta : r_\gamma \supset r_\alpha \circ r_\beta$.

Later, speaking about plurality of relations $\Re$ on $F$, which satisfies some of the above specified axioms, we will indicate numbers of the respective idioms for plurality symbol: the first will be the number of the axiom, which is responsible for reflexivity of the upper edge, the second – for symmetricity, and the third will indicate the axiom, which is responsible for transitivity of the upper edge. Thus, for instance, weakly reflexive, weakly symmetric and normally transitive plurality of relations is $\Re_{112}$-plurality.

Out of possible conditions for plurality elements, which ensure the property of the upper edge to be the relation of equivalency, it is necessary to identify those, which will ensure the possibility of constructive analysis of its properties, at preserving of initially built-in commonality in statement of the question. Let us start with the examples.

**4.4. Examples.** Hereinafter, $R, \overline{R}_+, R_+$ will be the designations for the set of real numbers and the subset of positive and strictly positive real numbers. Designations $Q, \overline{Q}_+, Q_+$ have analogous sense for subsets of rational numbers.

**1**. *Natural topology of real numbers:* let $F = R$, and $I = R_+$. We will consider that $y$ out of $F$ exists in relation $r_\alpha$ with $x$ out of $F$, if $|x - y| < \alpha$, where $\alpha \in I$.

Let us consider plurality of relations $\{r_\alpha\}_{\alpha \in I}$. It is clear that $x \in r_\alpha[x]$ at any $r_\alpha$, and each $r_\alpha$ is strongly symmetric. And finally, for any $\alpha, \beta$ one can find $\gamma > \alpha + \beta$ of such kind, that $r_\alpha \circ r_\beta \subset r_\gamma$; the opposite is also true. Hence, we got $\Re_{344}$-plurality of relations in $F$.

**2**. *Groups* (*of transformations*). Let in set $F$ some $\Re_{222}$-plurality of algebraic relations be preset. Hence, there exists $r_0 = id$; for any $r_\alpha$ one can find $r_\beta = r_\alpha^{-1}$; finally, for all $r_\alpha$ and $r_\beta$ one can find $r_\gamma = r_\alpha \circ r_\beta$. Recollecting associativity of composition of relations, we will get, that $\Re_{222}$-plurality of algebraic relations is a group (of transformations) ⌛

Analogous analysis of the other examples from topology and functional analysis show, that common conditions for elements of relations plurality, which induce miscellaneous mathematical structures, are: weak reflexivity, weak symmetric property, and normal transitivity. Moreover, *normal transitivity is guaranteed by definition of composition* for random relations. Later, the lower edge of plurality (in the first example it is trivial, and in the second it is primitive) is not an element of relations plurality. It produces the following:

**4.5. Definition.** $\Re_{112}$-plurality of relations in $F$ without minimal element is identified as binary relations filter.

When needed, we will portray the filter as pair $\Re_F = \{\Re_{112}; F\}$. Conformable to it, set $F_\Re$, in which the filter is defined, will be called filterable set and portrayed as pair $F_\Re = \{F; \Re_{112}\}$.

*Filter image* in point $x$ is the name for set $J_r[x]$. In different points $x$ out of domain of definition the filter can have different image (– equivalency classes, or, as we will say, – *strata*).

*Filter germ* $\Re_F = \{\Re_{112}; F\}$ will be the name for lower edge of relations plurality, which induces filter, with image $j_r[x]$ in some point $x \in F$.

*Sub-filter* of filter $\Re_F = \{\Re_{112}; F\}$ is the name for the filter, which is obtained by narrowing the domain of definition (relations plurality) of the initial filter: $\Re'_{F'} = \{\Re'_{112} \subset \Re_{112}; F' \subset F\}$.

When the filters under investigation have similar domain of definition, we will specify only the pluralities of relations corresponding them.

As a result of it, both investigation of the filters and all kinds of relations between them become the subject of the mathematics. Thus, variational calculus investigates mappings extremals of the set being filtered into natural topology filter of the real numbers; differential and integral equations investigate different kinds of mappings of different strata of functions sets, being filtered inside themselves, etc.

From this point on, much attention will be compelled to $\Re_{312}$-filters. It is just they, that will lead us to broader synthesis of the notion "topological space" and the other intriguing objects, which do not satisfy standard axioms of topological structures.